\documentclass[reqno,11pt]{amsart}

\usepackage{amsmath}
\usepackage{amsthm}
\usepackage{mathrsfs}
\usepackage{amsfonts}
\usepackage{parskip} 
\usepackage{hyperref}
\usepackage{epsfig}

\usepackage{amsmath,amssymb,amsfonts,amsthm,amstext,verbatim}
\usepackage{enumerate,color,graphicx,
}
\usepackage{psfrag}  
\usepackage{url}
\usepackage{mdwlist}   
\usepackage{mathrsfs}
\usepackage{hyperref}
\usepackage{epsfig}
\usepackage{etoolbox}

\makeatletter
\patchcmd{\@part}{\Huge}{\normalsize}{}{} 
\makeatother

\numberwithin{equation}{section}

\allowdisplaybreaks

\newtheorem{theorem}{Theorem}[section]

\newtheorem{proposition}[theorem]{Proposition}
\newtheorem{lemma}[theorem]{Lemma}

\newtheorem{definition}[theorem]{Definition}
\newtheorem{corollary}[theorem]{Corollary}

\newcommand{\ep}{\varepsilon}

\newcommand{\R}{\mathbb{R}}
\newcommand{\N}{\mathbb{N}}
\newcommand{\Z}{\mathbb{Z}}

\newcommand{\dist}{\vert\vert}

\newcommand{\Cgjthree}{C_4}

\newcommand{\eps}{\ep}

\newcommand{\psap}{\mathsf{P}}

\newcommand{\psaw}{\mathsf{W}}
\newcommand{\psawtd}{\mathsf{W}^{3,\mathsf{E}}}
\newcommand{\barp}{\hat{p}}
\newcommand{\carp}{\hat{c}}
\newcommand{\thethat}{\hat{\thet}}
\newcommand{\muhat}{\hat{\mu}}
\newcommand{\psaptd}{\mathsf{P}^{3,\mathsf{E}}}
\newcommand{\sawtd}{\mathrm{SAW}^{3,\mathsf{E}}}
\newcommand{\saptd}{\mathrm{SAP}^{3,\mathsf{E}}}

\newcommand{\esaw}{{\mathbb{E}_{\mathsf{W}}}}

\newcommand{\sapellward}{\mathrm{SAP}^{\mathfrak{l}}}

\newcommand{\saprward}{\mathrm{SAP}^{\mathfrak{r}}}

\newcommand{\reflect}{\mathcal{R}}

\newcommand{\lattv}{\Z^d}

\newcommand{\northeast}{\mathrm{NE}}
\newcommand{\eastnorth}{\mathrm{EN}}
\newcommand{\eastsouth}{\mathrm{ES}}
\newcommand{\westnorth}{\mathrm{WN}}
\newcommand{\southeast}{\mathrm{SE}}

\newcommand{\rightv}{\mathrm{Right}}

\theoremstyle{remark}

\newcounter{mycount}

\def\mik{1}
\newcommand\cpsfrag[2]{\ifnum\mik=1\psfrag{#1}{#2}\fi}

\newcommand{\saw}{\mathrm{SAW}}
\newcommand{\rhssaw}{\mathrm{RHSSAW}}

\newcommand{\sap}{\mathrm{SAP}}

\newcommand{\Cpcp}{\zeta}
\newcommand{\maceta}{\zeta}
\newcommand{\Cpoly}{C_1}
\newcommand{\Ch}{C_6}

\newcommand{\Cgjone}{C_2}
\newcommand{\Cgjtwo}{C_3}

\newcommand{\thet}{\theta}

\newcommand{\closes}{\text{ closes}}

\newcommand\Ga{\Gamma}     
\newcommand\ga{\gamma}     

\newcommand{\hidden}[1]{}

\newcommand\strongjoin{\mathsf{GlobalMJ}}
\newcommand\strongtdjoin{\mathsf{GlobalJoin}}
\newcommand\globaljoin{\mathsf{GJ}}
\newcommand\ymax{y_{{\rm max}}}
\newcommand\ymin{y_{{\rm min}}}
\newcommand\xmax{x_{{\rm max}}}
\newcommand\xmin{x_{{\rm min}}}

\newcommand\righttip{\eastsouth}

\newcommand{\modify}[1]{#1_{{\rm mod}}}

\newcommand{\highpolynum}{\mathsf{HPN}}

\newcommand{\proj}{\mathsf{Proj}}

\newcommand{\leftv}{\mathrm{Left}}
\newcommand{\lrpp}{\mathsf{LeftRightPolyPair}^{3,\mathsf{E}}}
\newcommand{\xspan}{x_{{\rm span}}}

\newcommand{\targetmany}{\mathsf{TMHA}}

\title[Self-avoiding polygons]{An upper bound on the number of \\ self-avoiding polygons via joining}
\date{}

\author[A.~Hammond]{Alan Hammond}
\address{Departments of Mathematics and Statistics, 
  U.C. Berkeley,
  Berkeley, CA, 94720-3840, U.S.A.}
\email{alanmh@stat.berkeley.edu}

\begin{document}
\maketitle

\vspace{-0.6cm}

\begin{abstract}
For $d \geq 2$ and $n \in \N$ even, let $p_n = p_n(d)$ denote the number of length~$n$ self-avoiding polygons in $\Z^d$ up to translation.  The polygon cardinality grows exponentially, and the growth rate  $\lim_{n \in 2\N} p_n^{1/n} \in (0,\infty)$ is called the connective constant and denoted by~$\mu$. Madras~[J. Statist. Phys. 78 (1995) no. 3--4, 681--699] has shown that $p_n \mu^{-n} \leq C n^{-1/2}$ in dimension $d=2$. Here we establish that
$p_n \mu^{-n} \leq n^{-3/2 + o(1)}$ for a set of even $n$ of full density when $d=2$.
We also consider a  certain variant of self-avoiding walk and argue that, when $d \geq 3$, an upper bound of  $n^{-2  + d^{-1} + o(1)}$  holds on a full density set for the counterpart in this variant model of this normalized polygon cardinality.  
\end{abstract}
\maketitle



\section{Introduction}

Self-avoiding walk was introduced in the 1940s by Flory and Orr \cite{Flory,Orr47} as a model of a long polymer chain in a system of such chains at very low concentration. It is well known among the basic models of discrete statistical mechanics for posing problems that are simple to state but difficult to solve. Two recent surveys are the lecture notes~\cite{BDGS11}  and \cite[Section 3]{Lawler13}. 

\subsection{The model} 

We will denote by $\N$ the set of non-negative integers.
Let $d \geq 2$. For $u \in \R^d$, let $\dist u \dist$ denote the Euclidean norm of $u$. 
A {\em walk} of length $n  \in \N$ with $n > 0$ is a map $\gamma:\{0,\cdots,n \} \to \lattv$ 
such that $\dist \gamma(i + 1) - \gamma(i) \dist = 1$ for each $i \in \{0,\cdots,n-1\}$.
An injective walk is called {\em self-avoiding}. 
A walk $\gamma \in \saw_n$ is said to close (and to be closing)
if $\dist \gamma(n) \dist = 1$. 
When the {\em missing edge} connecting $\gamma(n)$ and $\gamma(0)$
is added, a polygon results.
\begin{definition}
Let $\ga:\{0,\ldots,n-1\} \to \Z^d$ be a closing self-avoiding walk. For  $1 \leq i \leq n-1$, let $u_i$ denote the unordered nearest neighbour edge in $\Z^d$ with endpoints $\ga(i-1)$ and $\ga(i)$. Let $u_n$ denote $\ga$'s missing edge, with endpoints $\ga(n-1)$ and~$\ga(0)$. 
We call the collection of edges 
$\big\{ u_i: 1 \leq i \leq n \big\}$ the polygon of~$\ga$. A self-avoiding polygon in $\Z^d$ is defined to be any polygon of a closing self-avoiding walk in~$\Z^d$. The polygon's length is its cardinality. 
\end{definition}
We will usually omit the adjective self-avoiding in referring to walks and polygons. Recursive and algebraic structure has been used to analyse polygons in such domains as strips, as~\cite{BMB09} describes.

Note that the polygon of a closing walk has length that exceeds the walk's by one. Polygons have even length and closing walks, odd.


\subsection{Main results}

This article has two main conclusions, Theorems~\ref{t.polydev} and~\ref{t.threed}; both are upper bounds on the number of polygons of given length (the latter for a variant model). We present the two results now.

Define the polygon number $p_n$ to be the number of length~$n$ polygons up to translation. By (3.2.9) of~\cite{MS93}, 
$\lim_{n \in 2\N} p_n^{1/n} \in (0,\infty)$ exists; it called the connective constant and denoted by $\mu$. We
define the real-valued {\em polygon number deficit} exponent $\theta_n$ according to the formula
\begin{equation}\label{e.thetn}
 p_n = n^{-\thet_n} \cdot \mu^n \, \, \, \textrm{for $n \in 2\N$} \, .
\end{equation}

That $\liminf_{n \in 2\N} \thet_n$ is non-negative in any dimension~$d \geq 2$ will be reviewed in Section~\ref{s.thetalb}. 
Madras~\cite{Madras91} has proved a bound on the moment generating function of the sequence $\big\{ p_n: n \in 2\N \big\}$ which when $d=3$ would assert~$\lim_{n \in 2\N} \thet_n \geq 1$ were this limit known to exist.
More relevantly for us, he has shown in~\cite{Madras95} using a polygon joining technique that $\theta_n \geq 1/2 - o(1)$ for $d=2$. We develop this technique to prove a stronger lower bound valid for typical high~$n$.

\begin{definition}
The limit supremum density of 
a set $A$ of even, or odd, integers~is
$$
\limsup_n \frac{\big\vert A \cap [0,n] \big\vert}{\vert 2\N \cap [0,n] \vert} = \limsup_n \, n^{-1} \big\vert A \cap [0,2n] \big\vert \, .
$$
When the corresponding limit infimum density equals the limit supremum density, we naturally call it the limiting density. 
\end{definition}
\begin{theorem}\label{t.polydev}
 Let $d=2$. For any $\delta > 0$,
the  limiting density of the set of $n \in 2\N$ for which $\thet_n \geq 3/2 - \delta$ is equal to one. 
\end{theorem}

As we will discuss shortly, in the strongest form of this result consistent with predictions, the value $5/2$ would replace $3/2$.



 Madras and Slade \cite[Theorem 6.1.3]{MS93} have proved that $\thet_n \geq d/2 + 1$ in dimension $d \geq 5$ for {\em spread-out} models, in which the vertices of $\Z^d$ are connected by edges below some bounded distance.  There is a prospect that the method of proof of Theorem~\ref{t.polydev} may be applied in all dimensions $d \geq 2$, and the expected conclusion is that  $\thet_n \geq 2 - d^{-1} -  o(1)$ for a full density set of even $n \in \N$. Indeed, our second main conclusion, Theorem~\ref{t.threed}, is a result to this effect. 
 Certain technical difficulties arise as we try to apply our method in the higher dimensional case and, in order to largely circumvent these, we present  Theorem~\ref{t.threed} for a variant model.

\begin{definition}
The maximal edge local time of a nearest neighbour walk $\gamma:\{ 0,\cdots,n\} \to \Z^d$ is the maximum number of times that $\ga$ traverses an edge of $\Z^d$; more formally, it is the maximum cardinality of a subset $I \subseteq \{ 0,\cdots,n-1 \}$ such that the unordered sets  $\big\{ \ga(i),\ga(i+1) \big\}$ and $\big\{ \ga(j),\ga(j+1) \big\}$ coincide for each pair $(i,j) \in I^2$. Call $\ga$ $k$-edge self-avoiding if its maximal edge local time is at most~$k \in \N$.
Note that even $1$-edge self-avoiding walk satisfies a weaker avoidance constraint than does self-avoiding walk.

When considering (as we will) $3$-edge self-avoiding walks, we say that a walk $\ga$ as above closes if $\dist \ga(n) \dist = 0$.  Two such walks may be identified if they coincide after reparametrization by cyclic shift or reversal. A $3$-edge self-avoiding polygon is an equivalence class under this relation on closing walks. The length of such a polygon is the length of any of its members (and is $n$ if one of these members is $\ga$ as above).  Note that these definitions entail that
not only polygons but also
 closing walks have even length.
\end{definition}

For $n \in 2\N$, let $\barp_n$ 
denote the number of
$3$-edge self-avoiding polygons of length~$n$ up to translation. 
The connective constant $\lim_{n \in 2\N} \barp_n^{1/n}$ also exists for this model and we denote it by~$\muhat$. We define a real sequence $\big\{ \thethat_n : n \in \N \big\}$ so that
\begin{equation}\label{e.thetn.td}
 \barp_n = n^{-\thethat_n} \cdot \muhat^n \, \, \, \textrm{for $n \in 2\N$} \, .
\end{equation}
 
\begin{theorem}\label{t.threed}
Let $d \geq 2$. 
 For any $\delta > 0$,
the set of $n \in 2\N$ for which $\thethat_n \geq 2 - 1/d  - \delta$
has limiting  density equal to one. 
\end{theorem}

\subsection{Corollaries of the main results: the closing probability}

Let $\saw_n$ denote the set of self-avoiding walks~$\gamma$ of length~$n$ that start at $0$, i.e., with $\gamma(0) = 0$. 
We denote by $\psaw_n$ the uniform law on $\saw_n$. 
The walk under the law $\psaw_n$ will be denoted by $\Ga$. 
The {\em closing probability} is  $\psaw_n \big( \Gamma \closes \big)$.

Let the walk number $c_n$ equal the cardinality of $\saw_n$. By equations~(1.2.10) and~(3.2.9) of~\cite{MS93}, the limit 
$\lim_{n \in \N} c_n^{1/n}$ exists and coincides with~$\mu$, and we have~$c_n \geq \mu^n$.

The closing probability may be written in terms of the polygon and walk numbers.
 There are $2n$ closing walks whose polygon is a given polygon of length~$n$, since there are $n$ choices of missing edge and two of orientation. Thus,
\begin{equation}\label{e.closepc}
\psaw_n \big( \Gamma \closes \big) = \frac{2(n+1) p_{n+1}}{c_n} \, ,
\end{equation}
for any $n \in \N$ (but non-trivially only for odd values of $n$).
Since $c_n \geq \mu^n$,
Theorem~\ref{t.polydev} and~(\ref{e.closepc}) imply the next result.
\begin{corollary}\label{c.closingprob}
Let $d =  2$. For any $\eps > 0$, the set of $n \in 2\N + 1$ such that
$$
\psaw_n \big( \Gamma \closes \big) \leq n^{-1/2 + \eps} 
$$
has limiting density equal to one.
\end{corollary}

The walk (and polygon) cardinalities, and the closing probability, are notions that may be formulated for our variant model. Indeed, we may write $\carp_n$ and $\psawtd_n$ for the cardinality of, and uniform law on, the set $\sawtd_n$ of length~$n$ $3$-edge self-avoiding walks beginning at $0$.  
As we will explain in Section~\ref{s.threed}, we also have~$\carp_n \geq \muhat^n$. The counterpart 
$$
\psawtd_n \big( \Gamma \closes \big) = \frac{2 n \barp_n}{\carp_n} 
$$
to~(\ref{e.closepc}) is clearly valid
for any $n \in 2\N$. As such, the next inference is immediate from Theorem~\ref{t.threed}.
\begin{corollary}\label{c.threed}
Let $d \geq  2$. 
 For any $\delta > 0$,
the set of $n \in 2\N$ for which 
$$
\psawtd_n \big( \Gamma \closes \big) \leq n^{- (1 - 1/d) + \delta}
$$
has  limiting  density equal to one. 
\end{corollary}

\subsection{Further applications of the main results}

In~\cite{ontheprob}, an upper bound on the closing probability of $n^{-1/4 + o(1)}$ was proved in general dimension.  The method leading to that result was reworked in~\cite{snakemethodpattern} to prove the next result, which supersedes Corollary~\ref{c.closingprob}.

\begin{theorem}\label{t.closingprob}
Let $d \geq 2$. For any $\eps > 0$ and $n \in 2\N + 1$ sufficiently high,
$$
\psaw_n \big( \Gamma \closes \big) \leq n^{-1/2 + \eps} \, .
$$
\end{theorem}

The snake method is a probabilistic tool that is introduced in~\cite{snakemethodpattern} in order to obtain Theorem~\ref{t.closingprob}. (The article~\cite{snakemethodpattern}  is in large part a reworking of~\cite{ontheprob}, although this term for the method was not used in the original article~\cite{ontheprob}.) The snake method is applied via a technique of Gaussian pattern fluctuation in~\cite{snakemethodpattern}  to obtain Theorem~\ref{t.closingprob}. A second application of the snake method is made in~\cite{snakemethodpolygon}  in order to reach a stronger inference regarding the closing probability, valid when $d=2$, for a subsequence of even~$n$. In this case, the snake method is allied not with a pattern fluctuation technique but instead with the tool that is central to the proofs of the present article's principal results Theorems~\ref{t.polydev} and~\ref{t.threed}: the polygon joining technique, initiated by Madras in~\cite{Madras95}.   The conclusion that is reached in~\cite{snakemethodpolygon} 
 is now stated. 
\begin{theorem}\label{t.thetexist}
Let $d = 2$. 
\begin{enumerate} 
\item
For any $\eps > 0$, the bound 
$$
\psaw_n \big( \Gamma \closes \big) \leq n^{-4/7 + \eps} 
$$
holds on a set of $n \in 2\N + 1$ of limit supremum density at least $1/{1250}$. 
\item Define the {\em walk number excess} exponents $\xi_n$, $n \in \N$, via  $c_n = n^{\xi_n} \mu^n$, in counterpart to~(\ref{e.thetn}); they are non-negative by (1.2.10) of~\cite{MS93}.
Suppose that the limits $\thet : = \lim_{n \in 2\N} \thet_n$
and $\xi : = \lim_{n \in \N} \xi_n$
exist in $[0,\infty]$. 
Then $\thet + \xi \geq 5/3$. Since 
\begin{equation}\label{e.closepc.formula}
\psaw_n \big( \Gamma \closes \big) = n^{-\theta - \xi + 1 + o(1)} 
\end{equation}
as $n \to \infty$ through odd values of $n$
by (\ref{e.closepc}), the closing probability is seen to be bounded above by $n^{-2/3 + o(1)}$. 
\end{enumerate} 
\end{theorem} 
(When $\thet + \xi = \infty$, (\ref{e.closepc.formula}) should be interpreted as asserting a superpolynomial decay in $n$ for the left-hand side.)

The above result not only shares an important common element with Theorem~\ref{t.polydev} in the form of the technique used to prove the two results; the proof of Theorem~\ref{t.thetexist} also uses the statement of 
 Theorem~\ref{t.polydev}. 

It may be apparent that this paper shares with \cite{snakemethodpattern} and \cite{snakemethodpolygon}
certain combinatorial and probabilistic elements.
The present article, and the other two, may be read alone, but there may be also be value in viewing the results in unison. 
 The online article~\cite{CJC} presents the content of the three articles in a single work and includes some expository discussion. The online work connects the ideas; for example, its Section~3.3 elaborates the upcoming heuristic derivation of the polygon number deficit exponent lower bound~(\ref{e.thetonesided}) and uses this derivation as a mnemonic to discuss proofs of several of the main conclusions.

\subsection{Exponent prediction and hyperscaling relation}
 The limiting value $\thet = \lim_{n \in 2\N} \theta_n$ is predicted to exist and to satisfy a relation with the Flory exponent~$\nu$ for mean-squared radius of gyration. The latter exponent is specified by the putative formula 
\begin{equation}\label{e.nu}
       \esaw_n \, \dist \Gamma(n) \dist^2   = n^{2\nu + o(1)}  \, ,
\end{equation}
where $\esaw_n$ denotes the expectation associated with~$\psaw_n$ (and where note that $\Gamma(n)$ is the non-origin endpoint of $\Gamma$); in essence, $\dist \Ga(n) \dist$  is supposed to be typically of order $n^{\nu}$. 
The hyperscaling relation that is expected to hold between $\thet$ and $\nu$ is 
\begin{equation}\label{e.hyperscaling}
    \thet   =   1 +  d \nu
\end{equation}
where the dimension $d \geq 2$ is arbitrary. In $d = 2$, $\nu = 3/4$ and thus $\theta = 5/2$ is expected. That $\nu = 3/4$ was predicted by the Coulomb gas formalism \cite{Nie82,Nie84} and then by conformal field theory \cite{Dup89,Dup90}. 
Hara and Slade~\cite{HS92,HS92b} used the lace expansion to show that $\nu = 1/2$ when $d \geq 5$ by demonstrating that, for some constant $D \in (0,\infty)$,  $\esaw_n \, \dist \Gamma(n) \dist^2 - Dn$ is $O(n^{-1/4 + o(1)})$. This value of $\nu$ is anticipated in four dimensions as well, since $\esaw_n \, \dist \Gamma(n) \dist^2$ is expected to grow as $n \big( \log n \big)^{1/4}$. 
In fact, the continuous-time weakly self-avoiding walk in $d=4$ has been the subject of an extensive recent investigation of Bauerschmidt, Brydges and Slade.  In \cite{BBSlog}, a $\log^{1/4}$ correction to the susceptibility is derived,
relying on 
 rigorous renormalization group analysis developed in a five-paper series~\cite{BSRigRen1,BSRigRen2,BSRigRen3,BSRigRen4,BSRigRen5}.

\medskip

\noindent{\bf The structure of the paper.} 
The paper has four further sections.

The first sets up notation. The second, Section~\ref{s.thetalb}, gives an expository overview of the polygon joining technique and, we hope, provides a useful background for reading the later proofs; it is not logically needed later except for gathering some standard facts that will be used. Section~\ref{s.three} proves Theorem~\ref{t.polydev}. The argument is varied in Section~\ref{s.threed} in order to prove Theorem~\ref{t.threed}.

\medskip

\noindent{\bf Acknowledgments.} 
I am very grateful to a referee for a thorough discussion of the article~\cite{CJC}. Indeed, the present form of the two main theorems in the present article is possible on the basis of a suggestion made by this referee, (and this strengthened form has led to an improvement in Theorem~\ref{t.thetexist}(1), proved in~\cite{snakemethodpolygon}). 
I thank Hugo Duminil-Copin and Ioan Manolescu for many stimulating and valuable conversations about the central ideas in the paper. I thank Wenpin Tang for useful comments on a draft version. I would also like to thank Itai Benjamini for suggesting the problem of studying upper bounds on the closing probability.

This work was supported by  NSF grant DMS-$1512908$.


\section{Some general notation and tools}

\subsubsection{Multi-valued maps}\label{s.mvp}
For a finite set $B$, let $\mathcal{P}(B)$ denote its power set.
Let $A$ be another finite set. A {\em multi-valued map} from $A$ to $B$ is a function $\Psi: A \to \mathcal{P}(B)$. An arrow is a pair $(a,b) \in A \times B$ for which $b \in \Psi(a)$; such an arrow is said to be outgoing from $a$ and incoming to $b$. We consider multi-valued maps in order to find lower bounds on $\vert B \vert$, and for this, we need upper (and lower) bounds on the number of incoming (and outgoing) arrows. This next lemma is an example of such a lower bound.
\begin{lemma}\label{l.mvm}
Let $\Psi: A \to \mathcal{P}(B)$.
Set $m$ to be the minimum over $a \in A$ of the number of arrows outgoing from $a$, and $M$ to be the maximum over $b \in B$ of the number of arrows incoming to~$b$. Then
  $\vert B \vert \geq m M^{-1} \vert 
A \vert$.
\end{lemma}
\noindent{\bf Proof.} The quantities $M \vert B \vert$ and $m \vert A \vert$
are upper and lower bounds on the total number of arrows. \qed

\subsubsection{Denoting walk vertices and subpaths}
For $i,j \in \N$ with $i \leq j$, we write $[i,j]$ for  $\big\{ k \in \N:  i \leq k \leq j \big\}$. For a walk $\ga:[0,n] \to \Z^d$ and $j \in [0,n]$, we write $\ga_j$ in place of $\ga(j)$. For $0 \leq i \leq j \leq n$, $\ga_{[i,j]}$ denotes the subpath $\ga_{[i,j]}: [i,j] \to \Z^d$ given by restricting $\ga$. 

\subsubsection{Notation for certain corners of polygons}

The most important ideas in the article may be communicated by considering the two-dimensional case via the proof of Theorem~\ref{t.polydev}. The proof of Theorem~\ref{t.threed} is a variation. We now present notation specific to the two-dimensional case.

\begin{definition}\label{d.corners}
The Cartesian unit vectors are denoted by $e_1$ and $e_2$ and the coordinates of $u \in \Z^2$ by $x(u)$ and $y(u)$.
For a finite set of vertices $V \subseteq \Z^2$, we define the northeast vertex $\northeast(V)$ in $V$ to be that element of $V$ of maximal $e_2$-coordinate; should there be several such elements, we take $\northeast(V)$ to be the one of maximal $e_1$-coordinate. That is, $\northeast(V)$ is the uppermost element of $V$, and the rightmost among such uppermost elements if there are more than one. Using the four compass directions, we may similarly define eight elements of $V$, 
including the lexicographically minimal and maximal elements of $V$, $\mathsf{WS}(V)$ and $\mathsf{EN}(V)$. We extend the notation to any self-avoiding walk or polygon $\gamma$, writing for example $\northeast(\gamma)$ for $\northeast(V)$, where $V$ is the vertex set of $\gamma$. 
For a polygon or walk $\ga$, set $\ymax(\ga) = y\big(\northeast(\gamma)\big)$, $\ymin(\ga) = y\big(\southeast(\gamma)\big)$, $\xmax(\ga) = x\big(\eastnorth(\gamma)\big)$ and $\xmin(\ga) = x\big(\westnorth(\gamma)\big)$.
The height~$h(\ga)$ of $\ga$ 
 is $\ymax(\gamma) - \ymin(\gamma)$  and its width~$w(\gamma)$ is  $\xmax(\gamma) - \xmin(\gamma)$.
\end{definition}

\subsubsection{Polygons with northeast vertex at the origin}

For $n \in 2\N$, let $\sap_n$ denote the set of length $n$ polygons $\phi$ such that $\northeast(\phi) = 0$. The set $\sap_n$ is in bijection with equivalence classes of length~$n$ polygons where polygons are identified if one is a translate of the other. Thus, $p_n =   \vert \sap_n \vert$. 

We write $\psap_n$ for the uniform law on $\sap_n$. 
A polygon  sampled with law $\psap_n$ will be denoted by $\Ga$, as a walk with law $\psaw_n$ is.

There are $2n$ ways of tracing the vertex set of a polygon $\phi$ of length~$n$: $n$ choices of starting point and two of orientation. We now select one of these ways. Abusing notation, we may write $\phi$ as a map from $[0,n]$ to $\Z^2$, setting $\phi_0 = \northeast(\phi)$, $\phi_1 = \northeast(\phi) - e_1$, and successively defining $\phi_j$ to be the previously unselected vertex for which $\phi_{j-1}$ and $\phi_j$ form the vertices incident to an edge in $\phi$, with the final choice $\phi_n = \northeast(\phi)$ being made. Note that $\phi_{n-1} = \northeast(\phi) - e_2$.

\subsubsection{Cardinality of a finite set $A$} This is denoted by either $\# A$ or $\vert A \vert$.


\subsubsection{Plaquettes}

The shortest non-empty polygons contain four edges. Certain such polygons play an important role in several arguments and we introduce notation for them now.

\begin{figure}
    \begin{center}
      \includegraphics[width=0.75\textwidth]{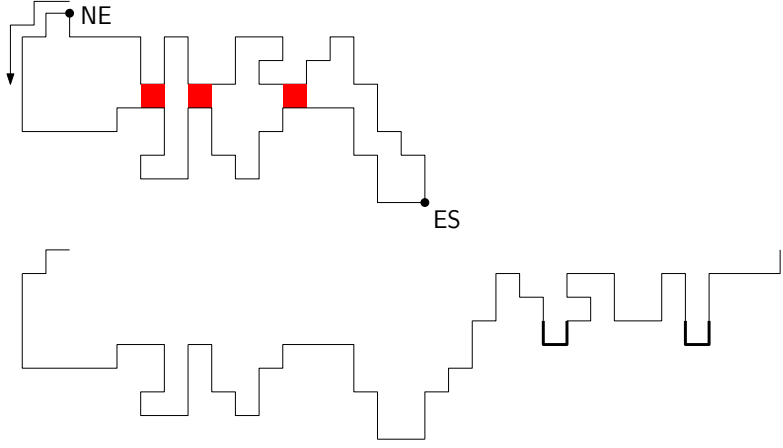}
    \end{center}
    \caption{In the upper sketch, a polygon $\phi$ with the beginning of the trajectory from $\northeast(\phi)$ indicated by an arrow. The polygon has three global join plaquettes, shaded red and labelled in the text $P^1$, $P^2$ and $P^3$ left to right (the order in which they are encountered counterclockwise along $\phi$ from $\northeast(\phi)$). In the lower picture, $\mathscr{S}_\kappa(\phi)$ with $\kappa = \big\{ 1,3 \big\}$ is depicted. The two bold subpaths indicate the reflections of the two diversions made by $\phi^{2,\kappa}$ that traverse $P^1$ and $P^3$ the long way around.}\label{f.reflect}
  \end{figure}

\begin{definition}
A plaquette is a polygon with four edges. Let $\phi$ be a polygon. A plaquette $P$ is called a join plaquette of~$\phi$ if $\phi$ and $P$ intersect at precisely the two horizontal edges of $P$. (The boundaries of the three shaded red squares of the polygon in the upcoming Figure~\ref{f.reflect} are join plaquettes of the polygon.)
Note that when $P$ is a join plaquette of $\phi$, the operation of removing the two horizontal edges in $P$ from $\phi$ and then adding in the two vertical edges in $P$ to~$\phi$ results in two disjoint polygons whose lengths sum to the length of $\phi$. We use symmetric difference notation and denote the output of this operation by $\phi \, \Delta \, P$.

The operation may also be applied in reverse: for two disjoint polygons $\phi^1$ and $\phi^2$, each of which contains one vertical edge of a plaquette~$P$, the outcome $\big( \phi^1 \cup \phi^2 \big) \, \Delta \, P$ of removing $P$'s vertical edges and adding in 
its horizontal ones is a polygon whose length is the sum of the lengths of $\phi^1$ and $\phi^2$.
\end{definition}

\section{Polygon number bounds via joining: an heuristic prelude}\label{s.thetalb}

In Section~2.1 of the text~\cite{MS93}, a presentation is made of the standard heuristic derivation of the relation~(\ref{e.hyperscaling}): note that equation~(1.4.14) of~\cite{MS93} is a representation of~(\ref{e.hyperscaling}) written in terms of the exponent $\alpha_{{\rm sing}} = 3 - \thet$.

When $d=2$, (\ref{e.hyperscaling}) asserts that $\theta = 1 + 2\nu$. A useful overview of our approach to proving Theorem~\ref{t.polydev}
is offered by giving in outline a heuristic derivation when $d=2$ of the weaker bound
\begin{equation}\label{e.thetonesided}
\thet \geq 1 +   \nu \, .
\end{equation}

To make the derivation, we hypothesise the existence of the limit
\begin{equation}\label{e.thetexist}
  \thet : = \lim_{n \in 2\N} \thet_n
\end{equation}
and suppose also that $\nu$ exists, given by ~(\ref{e.nu}).
We will derive~(\ref{e.thetonesided}) in three steps, arguing that in step $A$ that $\thet \geq 0$, in step $B$ that $\thet \geq \nu$ and concluding in step $C$. 

\medskip

\noindent{\bf Step A: $\theta \geq 0$.} By ~\cite[Theorem 3.2.3]{MS93}, we have that,
 for $n,m \in 2\N$, 
\begin{equation}\label{e.threepbound} 
 p_{n+m} \geq \tfrac{1}{d-1} p_n p_m \, .
\end{equation} 
Thus, the sequence $- \log \big( p_{2n}(d - 1)^{-1} \big)$, $n \in \N$, is subadditive, so that  Fekete's lemma \cite[Lemma 1.2.1]{St97} implies the existence of 
$\mu  = \lim_{n \in \N} p_{2n}^{1/{2n}}$ and the bound 
\begin{equation}\label{e.pbasicthree}
p_{2n} \leq (d-1)\mu^{2n} \, .
\end{equation}
 Thus, $\liminf_{n \in 2\N} \theta_n \geq 0$, completing step A. This step thus depends principally on the bound~(\ref{e.threepbound}), which is proved by a simple use of polygon joining. Our aim is to overview ideas for the upcoming proofs of our principal results, and  we consider only the case of $d=2$ and equal polygon length $m =n \in 2\N$ in explaining~(\ref{e.threepbound}).  
Consider a pair $\phi,\phi' \in \sap_n$. 
Relabel $\phi'$ by translating it so that $\westnorth(\phi')$ is one unit to the right of $\eastnorth(\phi)$. The plaquette $P$ whose upper left vertex is $\eastnorth(\phi)$ has one vertical edge in $\phi$ and one in $\phi'$. 
Note that $\chi : = \big( \phi \, \cup \, \phi' \big) \, \Delta \, P$ is a polygon of length~$2n$ which is either an element of~$\sap_{2n}$ or may be associated to such an element by making a translation. Moreover, the application $\big( \phi,\phi' \big) \to \chi$ is injective, because from~$\chi$ we can detect the location of the plaquette~$P$. The reader may wish to confirm this property, using that  $\phi$ and $\phi'$ have the same length. This injectivity implies that $p_{2n} \geq p_n^2$.

\medskip

\noindent{\bf Step B: $\theta \geq \nu$.} 
This step is an argument of Madras in~\cite{Madras95}.
When $d=2$, we will argue that 
\begin{equation}\label{e.madrasheur}
p_{2n} \geq n^{\nu} p_n^2 \, , 
\end{equation}
where recall that (\ref{e.nu}) specifies $\nu$. That $\theta \geq \nu$ follows directly (provided the two exponents exist).
The derivation of~(\ref{e.madrasheur}) develops the polygon joining argument in step A. The length $n$ polygons $\phi$ and $\phi'$ were joined in only one alignment, after displacement of $\phi'$ to a given location.
Madras argues under~(\ref{e.nu}) that there are at least an order of $n^\nu$ locations to which $\phi'$ may be translated and then attached to $\phi$. A total of $n^{\nu}$ distinct length $2n$ polygons results, and we obtain~(\ref{e.madrasheur}). 

 Where are these new locations for joining? Orient $\phi$ and $\phi'$ so that the height of each is at least its width; thus each height is at least of order $n^\nu$. Translate $\phi'$ vertically so that some vertices in $\phi$ and $\phi'$ share their $y$-coordinate, and then push  $\phi'$ to the right if need be so that this polygon is to the right of $\phi$. Then push $\phi'$ back to the left stopping just before the two polygons overlap. The two polygons contain vertices at distance one, and so it is plausible that in the locale of this vertex pair, we may typically find a plaquette $P$ whose left and right vertical edges are occupied by $\phi$ and $\phi'$. The operation in the proof of step A applied with plaquette $P$ then yields a length $2n$ polygon. 
This construction began with a vertical translation of $\phi'$,
with different choices of this translation resulting in different outcomes for the joined polygon; since there is an order of $n^\nu$ different heights that may be used for this translation, we see that the bound~(\ref{e.madrasheur}) results.

There is in fact a technical difficulty in implementing this argument: in some configurations, no such plaquette~$P$ exists. Madras developed a local joining procedure that overcomes this difficulty in two dimensions. His procedure will be reviewed shortly, in Section~\ref{s.madrasjoin}. 

\medskip

\noindent{\bf Step C: $\theta \geq 1 + \nu$.} 
To strengthen the conclusion to the form~(\ref{e.thetonesided}),
we argue heuristically in favour of a strengthening of~(\ref{e.madrasheur}),
\begin{equation}\label{e.joinstronger}
 p_{2n} \geq n^{\nu - o(1)} \sum_{j = - n/2}^{n/2} p_{n -j} p_{n+j}  \, .
\end{equation}

Expressed using the polygon number deficit exponents, we would then have $n^{-\thet_{2n}} \geq n^\nu \sum_{j = - n/2}^{n/2} (n-j)^{-\thet_{n -j}} (n+j)^{-\thet_j}$. Using (\ref{e.thetexist}), the bound~(\ref{e.thetonesided}) results.

To argue for (\ref{e.joinstronger}), note that, in deriving $p_{2n} \geq n^\nu p^2_n$, each polygon pair $(\phi,\phi') \in \sap_n \times \sap_n$ resulted in $n^\nu$ distinct length~$2n$ polygons. The length pair $(n,n)$ may be varied to be of the form $(n-j,n+j)$ for any $j \in [0,n/2]$. We are constructing a multi-valued map 
$$
\Psi: \bigcup_{\vert j \vert \leq n/2} \sap_{n-j} \times \sap_{n+j} \to \mathcal{P} \big( \sap_{2n} \big) \, 
$$
to the power set of $\sap_{2n}$ which associates to each polygon pair $(\phi,\phi')$ in the domain an order of $n^\nu$
elements of $\sap_{2n}$.  Were $\Psi$ injective, we would obtain~(\ref{e.joinstronger}). (The term {\em injective} is being misused: we mean that no two arrows of $\Psi$ are incoming to the same element of~$\sap_{2n}$.) The map is not injective but it is plausible that it only narrowly fails to be so: that is, abusing notation in a similar fashion, for typical $\chi \in \textrm{Range}(\phi)$, the cardinality of $\Psi^{-1}(\chi)$ is at most $n^{o(1)}$. A definition is convenient before we argue this.
\begin{definition}
Let $\phi$ be a polygon, and let $P$ be one of $\phi$'s join plaquettes. Let  $\phi^1$ and $\phi^2$ denote the disjoint polygons of which $\phi \, \Delta \, P$ is comprised.  If each of $\phi^1$ and $\phi^2$ has at least one quarter of the length of $\phi$, then we call $P$ a {\em macroscopic} join plaquette.
\end{definition}

To see that $\Psi$ is close to being injective, note that each pre-image of $\chi \in \textrm{Range}(\phi)$ under $\Psi$ corresponds to a macroscopic join plaquette of $\chi$.
Each macroscopic join plaquette entails a probabilistically costly macroscopic {\em four-arm} event, where four walks of length of order~$n$ must approach the plaquette without touching each other.
That $\chi$ belongs to $\textrm{Range}(\phi)$ amounts to saying that $\chi$ is an element of $\sap_{2n}$ having at least one macroscopic join plaquette. 
The four-arm costs 
make it plausible that a typical such polygon has only a few such plaquettes, gathered together in a small neighbourhood.
Thus, $\Psi$ is plausibly close to injective so that~(\ref{e.joinstronger}), and~(\ref{e.thetonesided}), results. 
(We will later be faced with arguing that $\Psi$ is indeed close to injective; although we will obtain a version of~(\ref{e.joinstronger}) via a form of near-injectivity, we will not obtain this near-injectivity by making rigorous the above, rather vague, heuristic.)

\section{Polygon joining when $d=2$, rigorously: proving Theorem~\ref{t.polydev}}\label{s.three}

In this section, we prove Theorem~\ref{t.polydev}, in which recall that~$d=2$. Our job is clear in light of the preceding heuristic derivation: we must make step $C$ rigorous.  
 The next proposition is a key step. It offers a rigorous interpretation of~(\ref{e.joinstronger}).

\begin{definition}\label{d.highpolynum}
For $\maceta > 0$, the set  $\highpolynum_\maceta \subseteq 2\N$ of $\maceta$-{\em high polygon number} indices is given by
$$
 \highpolynum_\maceta = \Big\{ n \in 2\N : p_n \geq n^{- \maceta} \mu^n \Big\} \, .
$$
\end{definition}



\begin{proposition}\label{p.polyjoin.newwork}
For any $\maceta > 0$, there is a constant $\Cpoly = \Cpoly(\maceta) > 0$ such that, for $n \in 2\N \cap \highpolynum_{\Cpcp}$,
\begin{equation}\label{e.polyjoin.newwork}
 p_n \geq \frac{1}{\Cpoly \log n} \sum_{j \in 2\N \, \cap \, [2^{i-1},2^i]} (n-j)^{1/2} p_j \, p_{n-j}  \, ,
\end{equation}
where $i \in \N$ is chosen so that $n \in 2\N \cap [2^i,2^{i+1}]$. 
\end{proposition}

In essence, the meaning of the proposition is that when  $n \in 2\N \cap \highpolynum_{\Cpcp}$, $p_n$ is at least a small constant multiple of $n^{1/2} (\log n)^{-1} \sum p_j p_{n-j}$, where the sum is over an interval of order $n$ indices $j$ around the value $n/2$ (though such a statement does not follow directly from the proposition).


Perhaps our proof of  Theorem~\ref{t.polydev} via Proposition~\ref{p.polyjoin.newwork}
can be refined to quantify the rate of convergence of the density one index set in the theorem.
However, the proposition in isolation is inadequate for proving the conclusion~$\thet_n \geq 3/2 - o(1)$ for all $n \in 2\N$, because this tool permits occasional spikes in the value of the $p_n$, as the sequence
$$
 \thet_n = \left\{
  \begin{array}{l l}
    1 & \quad \text{if $n$ is a power of $2$,}\\
    3/2 + \tfrac{1}{100} & \quad \text{if $n \in 2\N$ is otherwise,}
  \end{array} \right.
$$
demonstrates.

This section has five subsections. The first four present tools needed for the proof of Proposition~\ref{p.polyjoin.newwork}, with the proof of this result in the fourth. The fifth proves Theorem~\ref{t.polydev} as a consequence of the proposition.

Proving Proposition~\ref{p.polyjoin.newwork} amounts to making step C rigorous. We need to define in precise terms the polygon joining technique that we will use to do this, and, in  
 the first subsection, we specify the local details of the technique due to Madras that we will use. The heuristic argument in favour of~(\ref{e.joinstronger})
depended on the near-injectivity of the multi-valued map~$\Psi$, which was argued by making a case for the sparsity of macroscopic join plaquettes.
In the second subsection, we present
Proposition~\ref{p.globaljoin} and Corollary~\ref{c.globaljoin},
our rigorous versions of this sparsity claim. In the third subsection, we set up the apparatus needed to specify our joining mechanism~$\Psi$, and in the fourth, we define and analyse the mechanism (and so obtain Proposition~\ref{p.polyjoin.newwork}).

\subsection{Madras' polygon joining procedure}\label{s.madrasjoin}

When a pair of polygons is close, there may not be a plaquette whose vertical edges are divided between the two elements of the pair.  We now recall Madras' joining technique which works in a general way for such pairs. 

\medskip

Consider two polygons $\tau$ and $\sigma$ of lengths $n$ and $m$ for which the intervals 
\begin{equation}\label{e.intint}
\textrm{$\big[ \ymin(\tau) - 1 , \ymax(\tau) + 1 \big]$ and  $\big[ \ymin(\sigma) - 1 ,\ymax(\sigma) + 1 \big]$ intersect} \, .
\end{equation}
 Madras' procedure joins  $\tau$ and $\sigma$ to form a new polygon of length $n + m + 16$ in the following manner. 

First translate $\sigma$ to the right by far enough that the $x$-coordinates of the vertices of this translate are all strictly greater than all of those of $\tau$. Now shift $\sigma$ to the left step by step until the first time at which there is a pair of vertices, one in $\tau$ and the other in the $\sigma$-translate, that share an $x$-coordinate and whose $y$-coordinates differ by at most two; such a moment necessarily occurs, by the assumption~(\ref{e.intint}). Write $\sigma' = \sigma + T_1 e_1$ (with $T_1 \in \Z$) for this particular horizontal translate of $\sigma$. There is at least one vertex $z \in \Z^2$ such that the set $\big\{ z - e_2,z, z + e_2 \big\}$ contains a vertex of $\tau$ and a vertex of $\sigma'$. The set of such vertices contains at most one vertex with any given $y$-coordinate. Denote by $Y$ the vertex $z$ with the maximal $y$-coordinate. 

\begin{figure}
    \begin{center}
      \includegraphics[width=0.6\textwidth]{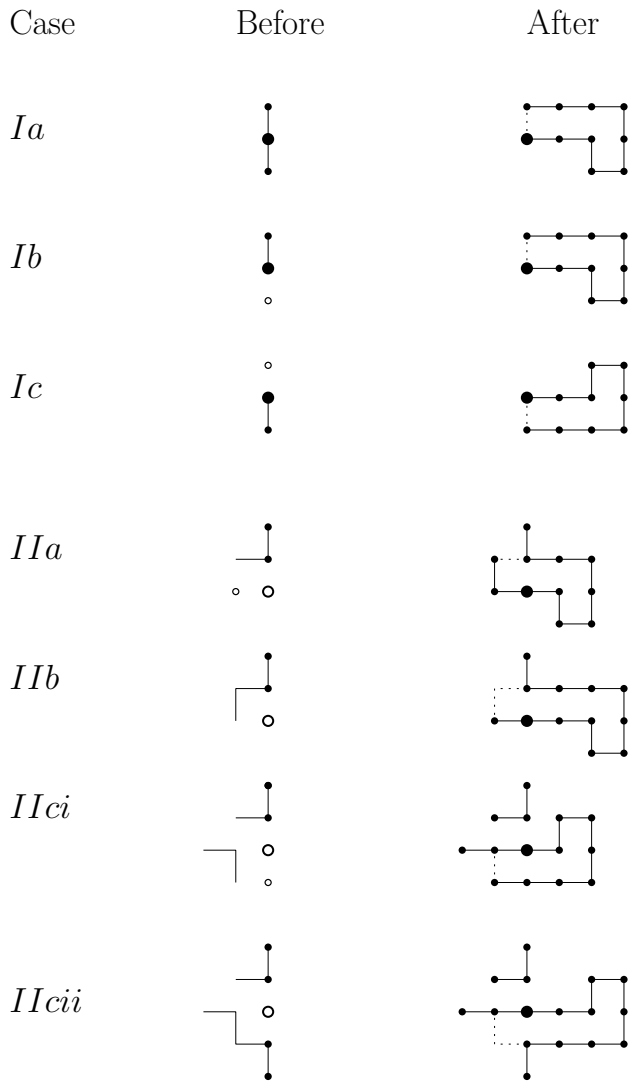}
    \end{center}
    \caption{Changes made near the vertex $Y$ in a polygon $\tau$ to produce $\modify{\tau}$ are depicted. The second column depicts $\tau$ around~$Y$, with $Y$ indicated by a large disk or circle; a disk denotes a vertex that belongs to $\tau$ and a circle one that does not; the line segments denote edges of $\tau$. In the third column, the modified polygon $\modify{\tau}$ is shown in the locale of $Y$. The vertex $Y$ is shown as a black disk. Black line segments are edges in $\tau$ or $\modify{\tau}$ and dotted line segments on the right are edges in $\tau$ that are removed in the formation of~$\modify{\tau}$. Several cases are not depicted. These may be labelled cases IIIa, IIIb, IIIci and IIIcii. In each case, the picture of $\tau$ and $\modify{\tau}$ is formed by reflecting the counterpart case II picture horizontally through~$Y$.}\label{f.madrascases}
  \end{figure}

Madras now defines a modified polygon $\modify{\tau}$, which is formed from $\tau$ by changing its structure in a neighbourhood of $Y \in \Z^2$. Depending on the structure of $\tau$ near~$Y$, either two edges are removed and ten edges added to form $\modify{\tau}$ from $\tau$, or one edge is removed and nine are added. As such, $\modify{\tau}$ has length $n + 8$. The rule that specifies $\modify{\tau}$ is recalled from~\cite{Madras95} in Figure~\ref{f.madrascases}. 

A modified polygon formed from $\sigma'$ is also defined. Rotate $\sigma'$ about the vertex~$Y$ by $\pi$ radians to form a new polygon $\sigma''$. Form $\modify{\sigma''}$ according to the same rules, recalled in Figure~\ref{f.madrascases}. Then rotate back the outcome by $\pi$ radians about $Y$ to produce the modification of 
$\sigma'$, which to simplify notation we denote by $\modify{\sigma}$. 

Writing $Y = (Y_1,Y_2)$ in place of $\big( x(Y) , y(Y) \big)$, note that no vertex of $\tau$ belongs to the {\em right corridor}
$\big\{ Y_1 + 1 , Y_1 + 2, \cdots \big\} \times  \big\{ Y_2 - 1 , Y_2 , Y_2 + 1 \big\}$, the region
that lies strictly to the right of $\big\{ Y-e_2,Y,Y+e_2 \big\}$. (Indeed, it is this fact that implies that $\modify{\tau}$ is a polygon.)
Equally, no vertex of $\sigma'$ belongs to the {\em left corridor}  
$\big\{ \cdots , Y_1 - 2 , Y_1 -1 \big\} \times \big\{ Y_2 - 1 , Y_2 , Y_2 + 1 \big\}$.

Note that the polygon $\modify{\tau}$ extends $\tau$ to the right of $Y$ by either two or three units inside the right corridor (by two in case IIa, IIci or IIIci and by three otherwise). Likewise $\modify{\sigma}$ extends $\sigma'$ to the left of $Y$ by either two or three units in the left corridor (by two when $\sigma''$ satisfies case IIa, IIci or IIIci and by three otherwise).

Note from Figure~\ref{f.madrascases} that, in each case, $\modify{\tau}$ contains two vertical edges that cross the right corridor at the maximal $x$-coordinate adopted by vertices in $\modify{\tau}$ that lie in this corridor. Likewise, $\modify{\sigma}$ contains two vertical edges that cross the left corridor at the minimal $x$-coordinate adopted by vertices in $\modify{\sigma}$ that lie in the left corridor.

Translate $\modify{\sigma}$ to the right by $T_2$ units, where $T_2$ equals
\begin{itemize}
\item five when one of cases IIa, IIci and IIIci obtains for both $\tau$ and $\sigma''$;
\item six  when one of these cases obtains for exactly one of these polygons;
\item seven when none of these cases holds for either polygon. 
\end{itemize}
Note that $\modify\tau$ and $\modify\sigma + T_2 e_1$ are disjoint polygons such that, for some pair of vertically adjacent plaquettes $(P^1,P^2)$ in the right corridor (whose left sides have $x$-coordinate either $Y_1 + 2$ or $Y_1 + 3$), the edge-set of $\modify\tau$ intersects the plaquette pair on the two left sides of $P^1$ and $P^2$, while the edge-set of $\modify\sigma + T_2 e_1$ intersects this pair on the two right sides of $P^1$ and $P^2$. Let $P^1$ denote the upper element of this plaquette pair.

The polygon that Madras specifies as the join of $\tau$ and $\sigma$ is given by $\big( \modify{\tau} \cup (\modify\sigma + T_2 e_1) \big) \, \Delta \, P^1$. 
Note that, to form the join polygon, $\sigma$ is first horizontally translated by $T_1$ units to form $\sigma'$, modified locally to form $\modify\sigma$, and then further horizontally translated by $T_2$ units to produce the polygon $\modify\sigma + T_2 e_1$ that is joined onto $\modify\tau$. Thus, $\sigma$ undergoes a horizontal shift by $T_1 + T_2$ units as well as a local modification before being joined with~$\modify\tau$.

  \begin{figure}
    \begin{center}
      \includegraphics[width=0.75\textwidth]{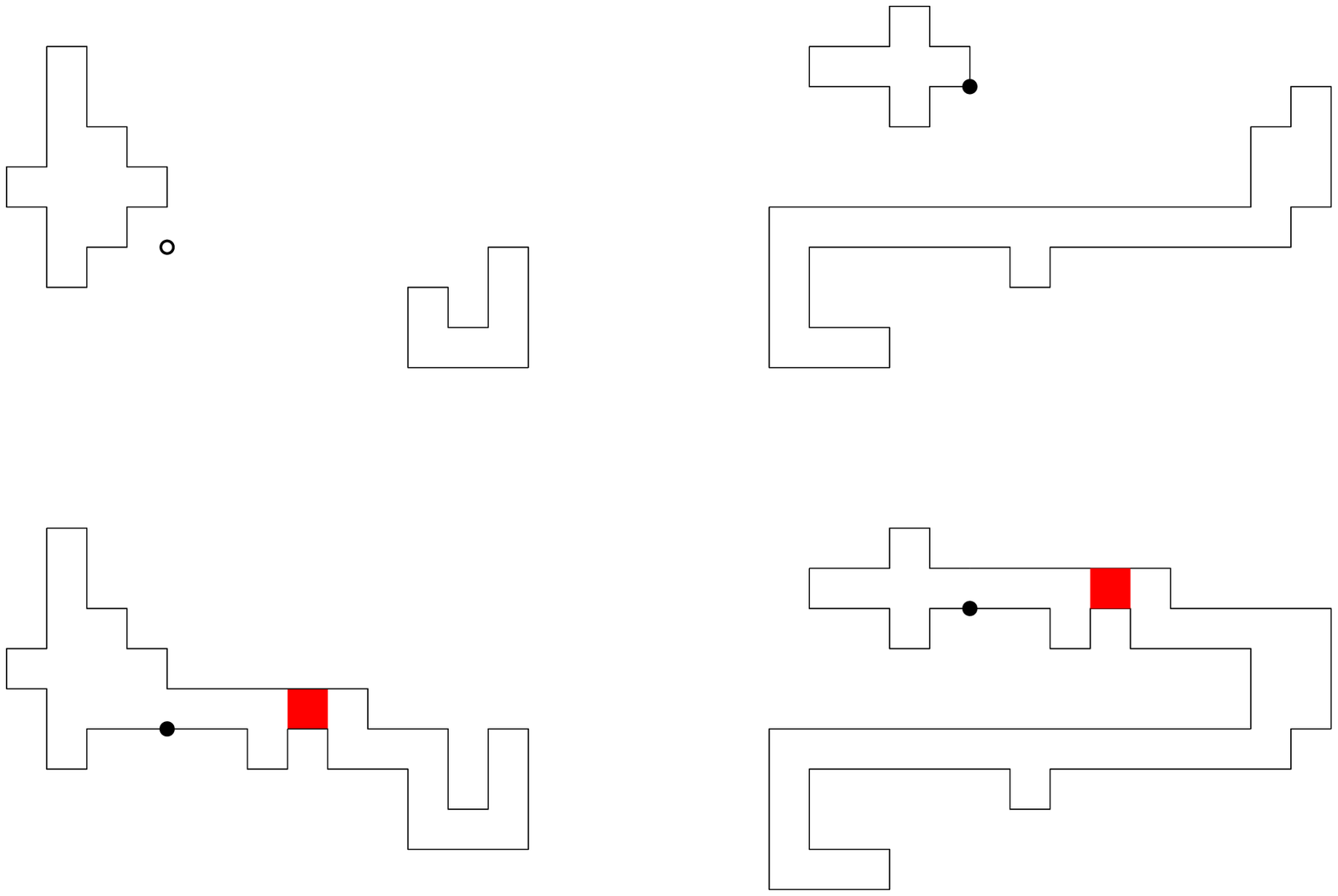}
    \end{center}
    \caption{Two pairs of Madras joinable polygons and the Madras join polygon of each pair. The vertex $Y$ is indicated by a circle or a disk in each picture. The junction plaquette is shaded red. On the left, $\tau$ satisfies case IIb and $\sigma''$, case IIa; on the right, Ib and~IIb.}\label{f.madrasjoinable}
  \end{figure}

\begin{definition}
For two polygons $\tau \in \sap_n$ and $\sigma \in \sap_m$ satisfying~(\ref{e.intint}), define the {\em Madras join polygon} 
$$
J(\tau,\sigma) = 
\Big( \modify{\tau} \cup \big( \modify\sigma + T_2 \, e_1 \big) \Big) \, \Delta \, P^1 \, \in \,  \sap_{n + m + 16} \, .
$$

The plaquette $P^1$ will be called the {\em junction} plaquette.

Such polygons $\tau$ and $\sigma$ are called {\em Madras joinable}  
if $T_1 + T_2 = 0$: that is, no horizontal shift is needed so that $\sigma$ may be joined to $\tau$ by the above procedure. Note that the modification made is local in this case: $J(\tau,\sigma) \Delta \big( \tau \cup \sigma \big)$ contains at most twenty edges. See Figure~\ref{f.madrasjoinable}.
\end{definition}

\subsection{Global join plaquettes are few}

Recall that in step C of the derivation of~(\ref{e.thetonesided}), the near-injectivity of the multi-valued map $\Psi$ was argued as a consequence of the sparsity of macroscopic join plaquettes. We now present in~Corollary~\ref{c.globaljoin} a rigorous counterpart to this sparsity assertion. In the rigorous approach, we use a slightly different definition to the notion of macroscopic join plaquette. 


\begin{definition}\label{d.globaljoin}
For $n \in 2\N$, let $\phi \in \sap_n$. 
A join plaquette $P$ of $\phi$ is called {\em global} if the two polygons 
comprising $\phi \, \Delta \, P$ may be labelled $\phi^\ell$ and $\phi^r$ in such a way that 
\begin{itemize}
\item every rightmost vertex in $\phi$ is a vertex of $\phi^r$;
\item and  $\northeast(\phi)$ is a vertex of $\phi^\ell$.
\end{itemize}
Write $\globaljoin_\phi$ for the set of global join plaquettes of the 
polygon~$\phi$.
\end{definition}

\begin{proposition}\label{p.globaljoin}
There exists $c > 0$ such that, for $n \in 2\N$ and any $k \in \N$,  
$$
  \# \Big\{ \phi \in \sap_{n} : \big\vert \globaljoin_\phi \big\vert \geq k \Big\} \leq  c^{-1} 2^{-k \mu^{-2}/2}    \mu^n   \, .
$$
\end{proposition}
\begin{corollary}\label{c.globaljoin}
For all $\maceta > 0$ and $\Cgjone > 0$, there exist $\Cgjtwo, \Cgjthree > 0$ such that, for $n \in \highpolynum_\maceta$,
$$
\psap_{n} \Big( \, \big\vert \globaljoin_\Gamma \big\vert \geq \Cgjtwo \log n \Big)
 \leq \Cgjthree n^{- \Cgjone} \, .
$$
\end{corollary}

The next lemma will be used in the proof of Proposition~\ref{p.globaljoin}.
\begin{lemma}\label{l.globaloneedge}
Let $n \in 2\N$ and $\phi \in \sap_n$.  Writing $j \in [0,n]$ so that $\eastsouth(\phi) = \phi_j$, consider the two subpaths
$\phi_{[0,j]}$ and $\phi_{[j,n]}$, the first starting at $\northeast(\phi) = 0$ and the second ending there. Each of these paths contains precisely one of the two horizontal edges of any element in $\globaljoin_\phi$.
\end{lemma} 
\noindent{\bf Proof.} For given $\phi \in \sap_n$, let $P \in \globaljoin_\phi$. We may decompose $\phi \, \Delta \, P$ as $\phi^\ell \cup \phi^r$ in accordance with Definition~\ref{d.globaljoin}. We then have that $\northeast(\phi)$ is a vertex of $\phi^{\ell}$ and $\eastsouth(\phi)$ a vertex of $\phi^{r}$. The path $[0,n] \to \Z: j \to \phi_j$ leaves $0$ to follow  $\phi^{\ell}$ until passing through an edge in $P$ to arrive in~$\phi^{r}$, tracing this polygon until passing back through the other horizontal edge of $P$ and following $\phi^{\ell}$ until returning to $\northeast(\phi)$. It is during the trajectory in $\phi^{r}$ that the visit to $\eastsouth(\phi)$ is made. \qed    

\medskip

\noindent{\bf Proof of Proposition~\ref{p.globaljoin}.}
For $n \in \N$, an element $\gamma \in \saw_n$
is called a {\em half-space} walk if $y(\gamma_i) \leq 0$ for each $i \in [0,n]$. We call a half-space walk $\ga$ {\em returning} if,
after the last visit that $\ga$ makes to the lowest $y$-coordinate that this walk attains, $\ga$ makes a unique visit to the $x$-axis, with this occurring at its endpoint $\ga_n$.
Let $\rhssaw_n \subset \saw_n$ denote the set of length~$n$ returning half-space walks. We will first argue that, for any $n \in \N$, 
\begin{equation}\label{e.rhswalkbound}
\vert \rhssaw_n \vert \leq \mu^n \, .   
\end{equation}
To see this, we consider a map $\rhssaw_n \to \saw_n$. Let $\ga \in \rhssaw_n$, and let $j \in [0,n]$
be the index of the final visit of $\ga$ to the lowest $y$-coordinate visited by~$\ga$. The image of $\ga$ is defined to be the concatenation of $\ga_{[0,j]}$ and the reflection of $\ga_{[j,n]}$ through the horizontal line that contains $\ga_j$.
(We have not defined concatenation but hope that the meaning is evident.) Our map is  injective because, given an element in its image, the horizontal coordinate of the line used to construct the image walk may be read off from that walk; (the coordinate is one-half of the $y$-coordinate of the image walk's non-origin endpoint). Its image lies in the set of {\em bridges} of length~$n$, where a bridge is a walk whose starting point has maximal $y$-coordinate and whose endpoint uniquely attains the minimal $y$-coordinate. 
The set of bridges of length~$n$ beginning at the origin 
has cardinality at most $\mu^n$: this classical fact, which follows from (1.2.17) in~\cite{MS93} by symmetries of the Euclidean lattice, is proved by a superadditivity argument with similarities to the proof of~(\ref{e.threepbound}). Thus, considering this map proves~(\ref{e.rhswalkbound}).

Noting these things allow us to reduce the proof of the proposition to verifying the following assertion. 
There exists $c > 0$ such that, for $\delta \in (0,1/2)$, $n \in 2\N$ and any $k \in \N$,  
\begin{equation}\label{e.globaljoin.reduce}
    \#  \, \rhssaw_{n + 2\lfloor \delta k \rfloor}  \, \geq \, c \, \delta^{-\delta k} \cdot \# \Big\{ \phi \in \sap_{n} : \big\vert \globaljoin_\phi \big\vert \geq k \Big\}  \, .
\end{equation}
Indeed, applying (\ref{e.rhswalkbound}) with the role of $n$ played by $n + 2\lfloor \delta k \rfloor$, we see from~(\ref{e.globaljoin.reduce}) that
$$
\# \Big\{ \phi \in \sap_{n} : \big\vert \globaljoin_\phi \big\vert \geq k \Big\} \leq c^{-1} \big( \mu^2 \delta \big)^{\delta k} \mu^n \, .
$$
Setting $\delta = \mu^{-2}/2$, we obtain Proposition~\ref{p.globaljoin}.

To complete the proof of the proposition, we must prove~(\ref{e.globaljoin.reduce}), and this we now do.
Let $\phi \in \sap_n$. Setting $j \in [0,n]$ so that 
$\eastsouth(\phi) = \phi_j$ (as we did in Lemma~\ref{l.globaloneedge}),
write $\phi^1 =  \phi_{[0,j]}$ and $\phi^2 = \phi_{[j,n]}$. 
Writing $\reflect^2_z$ for reflection in the vertical ($e_2$-directed) line that passes through $z \in \Z^2$,
define a map $\mathscr{S}: \sap_n \to \rhssaw_n$ to be the concatenation
$$
\mathscr{S}(\phi) = \phi^1 \circ \reflect^2_{\eastsouth(\phi)}(\phi^2) \, .
$$
By Lemma~\ref{l.globaloneedge}, each of $\phi^1$ and $\phi^2$ traverses precisely one horizontal edge of each of $\phi$'s global join plaquettes. Set $r = \# \globaljoin_\phi$ and enumerate $\globaljoin_\phi$ by the sequence $\big( P^1,\cdots,P^r \big)$ (in an arbitrary order; for example, in the order in which $\phi^1$ traverses an edge of each plaquette). 
For each $j \in \{1,\cdots,r \}$, let $(s_j,f_j)$ denote the unique edge in $P^j$ traversed by $\phi^2$. 
Consider the path formed by modifying $\phi^2$ so that the one-step subpath $(s_j,f_j)$ is replaced by a three-step subpath from $s_j$ to $f_j$ that traverses the plaquette $P^j$ using its three edges other than $(s_j,f_j)$. The modification may be made iteratively for several choices of $j \in \{ 1,\cdots,r \}$, and the outcome is independent of the order in which the modifications are made. In this way, we may define a modified path $\phi^{2,\kappa}$ for each $\kappa \subseteq \{1,\cdots,r\}$, under which the modified route is taken along plaquettes $P^j$ precisely when $j \in \kappa$. Note that $\phi^{2,\kappa}$ is a self-avoiding walk whose length exceeds $\phi^2$'s by $2\vert \kappa \vert$: it is self-avoiding because this walk differs from $\phi^2$
by several disjoint replacements of one-step subpaths by three-step alternatives, and, in each case, the two new vertices visited in the alternative route are vertices in $\phi^1$, and, as such, cannot be vertices in $\phi^2$.

Note further that the intersection of the edge-sets of $\phi^1$ and $\phi^{2,\kappa}$ equals $\cup_{j \in \kappa} E(\phi^1) \cap E(P^j)$ (where the sets in the union are each singletons).

For each $\kappa \subseteq \{ 1,\cdots,r\}$, define $\mathscr{S}_\kappa(\phi) \in \rhssaw_{n + 2 \vert \kappa \vert}$,
$$
\mathscr{S}_\kappa(\phi) = \phi^1 \circ \reflect^2_{\eastsouth(\phi)}(\phi^{2,\kappa}) \, .
$$

Recall that $\delta \in (0,1/2)$ and that $k \in \N$ is given. Consider the multi-valued map 
$$
\Psi = \Psi_k : \Big\{ \phi \in \sap_n: \vert \globaljoin_\phi \vert \geq k \Big\}
 \to \mathcal{P} \big( \rhssaw_{n + 2 \lfloor \delta k \rfloor} \big)
$$
that associates to each $\phi \in \sap_n$ with $\vert \globaljoin_\phi \vert \geq k$ the set 
$$
\Psi(\phi) = 
\Big\{ \mathscr{S}_\kappa(\phi): \kappa \subseteq \globaljoin_\phi \, , \,  \vert \kappa \vert = \lfloor \delta k \rfloor \Big\} \, ,
$$ 
where here we abuse notation and identify a subset of $\globaljoin_\phi$ with its set of indices under the given enumeration of $\globaljoin_\phi$.

Note that, for some constant $c > 0$ that is independent of $\delta \in (0,1/2)$ and for all $k \in \N$, 
$$
\big\vert \Psi(\phi) \big\vert =  {\vert \globaljoin_\phi \vert \choose \lfloor \delta k \rfloor} \geq  {k \choose \lfloor \delta k \rfloor} \geq c \, \delta^{-\delta k} \, .
$$
(In the latter displayed inequality, the factor of $k^{-1/2}$
that appears via Stirling's formula has been cancelled against other omitted terms. This detail is inconsequential and the argument is omitted.)

Note that, for any $\ga \in \saw_{n + 2 \lfloor \delta k \rfloor}$, the preimage $\Psi^{-1} (\gamma)$ is either the empty-set or a singleton. Indeed, if $\ga \in \Psi(\phi)$ for some $\phi \in \sap_n$ with $\vert \globaljoin_\phi \vert \geq k$, then $\phi$ may be recovered from $\ga$ as follows:
\begin{itemize}
\item the coordinate $x\big(\eastsouth(\phi)\big)$ equals $x\big( \ga_{n + 2 \lfloor \delta k \rfloor} \big)/2$;
\item the vertex $\eastsouth(\phi)$ is the lowest among the vertices of $\ga$ having the above $x$-coordinate;
\item setting $j \in [0,n]$ so that $\ga_j$ is this vertex, consider the non-self-avoiding walk $\ga_{[0,j]} \circ \reflect_{\ga_j}^2\big( \ga_{[j,n]} \big)$. This walk begins and ends at $0$. There are exactly $\lfloor \delta k \rfloor$ instances where the walk traverses an edge twice. In each case, the three-step journey that the walk makes in the steps preceding, during and following the second crossing of the edge follow three edges of a plaquette. Replace this journey by the one-step journey across the remaining edge of the plaquette, in each instance. The result is~$\phi$.
\end{itemize}

We may thus use the multi-valued map Lemma~\ref{l.mvm} to find a lower bound on the cardinality of $\rhssaw_{n + 2 \lfloor \delta k \rfloor}$. The resulting bound is precisely~(\ref{e.globaljoin.reduce}). This completes the derivation of Proposition~\ref{p.globaljoin}.   \qed

\medskip

\noindent{\bf Proof of Corollary~\ref{c.globaljoin}.}
Using
\begin{equation}\label{e.globalk}
\psap_n \Big(  \vert \globaljoin_\phi \vert \geq k \Big) = p_n^{-1} \cdot \# \Big\{ \phi \in \sap_n: \vert \globaljoin_\phi \vert \geq k \Big\} \, ,
\end{equation}
and Proposition~\ref{p.globaljoin},
we find that if $n \in \highpolynum_{\maceta}$ then 
$$
\psap_n \big(  \vert \globaljoin_\phi \vert \geq k \big) \leq c^{-1} \, 2^{- k \mu^{-2}/2}  n^{\Cpcp} 
$$
for $k \in \N$.
Set $\Cgjtwo = 2 \mu^2 \big( \log 2 \big)^{-1} \big( \Cgjone + \Cpcp + 1 \big)$. Then $n \geq c^{-1}$ implies that, for $k \geq \Cgjtwo \log n$,
$$
\psap_n  \big(  \vert \globaljoin_\phi \vert \geq k \big) \leq   c^{-1} n^{-\Cgjone} \, .
$$
Setting $\Cgjthree = c^{-1}$, we obtain Corollary~\ref{c.globaljoin}. 
 \qed

\subsection{Preparing for joining surgery: left and right polygons}\label{s.refpert.newwork}

When we prove Proposition~\ref{p.polyjoin.newwork} in the next subsection, polygons from a certain set $\saprward_{n-j}$
will be joined to others in another set $\sapellward_j$.
We think of the former as being joined to the latter on the right, so that the superscripts indicate a handedness associated to the joining.

Anyway, in this subsection, we specify these polygon sets, give a lower bound on their size in Lemma~\ref{l.polysetbound.newwork}, and then explain in Lemma~\ref{l.strongjoin} 
how there are plentiful opportunities for joining pairs of such polygons in a surgically useful way (so that the join may be detected by virtue of its global nature).

Let $\phi$ be a polygon. Recall from Definition~\ref{d.corners} the notation $\ymax(\phi)$ and $\ymin(\phi)$, as well as the height~$h(\phi)$ and width~$w(\phi)$.

\begin{definition}\label{d.leftright}
For $n \in 2\N$, let $\sapellward_n$ denote the set of {\em left} polygons $\phi \in \sap_n$
such that
\begin{itemize}
\item $h(\phi) \geq w(\phi)$ (and thus, by a trivial argument, $h(\phi) \geq n^{1/2}$),
\item and $y\big(\righttip(\phi) \big) \leq \tfrac{1}{2} \big( \ymin(\phi) + \ymax(\phi) \big)$. 
\end{itemize}
Let $\saprward_n$ denote the set of {\em right}  polygons $\phi \in \sap_n$
such that
\begin{itemize}
\item $h(\phi) \geq w(\phi)$.
\end{itemize}
\end{definition}
\begin{lemma}\label{l.polysetbound.newwork}
For $n \in 2\N$,
$$
\big\vert \sapellward_n \big\vert \geq \tfrac{1}{4} \cdot \big\vert \sap_n \big\vert \, \, \, 
\textrm{and}  \, \, \,
\big\vert \saprward_n \big\vert \geq  \tfrac{1}{2} \cdot \big\vert \sap_n \big\vert \, .
$$
\end{lemma}
\noindent{\bf Proof.} An  element $\phi \in \sap_n$
not in $\saprward_n$ is brought into this set
by right-angled rotation. 
If, after the possible rotation, it is not in $\sapellward_n$, it may 
brought there by reflection in the $x$-axis. \qed

\begin{definition}\label{d.globallymj}
A Madras joinable polygon pair $(\phi^1,\phi^2)$
is called globally Madras joinable if the junction plaquette of the join polygon $J(\phi^1,\phi^2)$ is a global join plaquette of  $J(\phi^1,\phi^2)$.
\end{definition}
Both polygon pairs in the upper part of Figure~\ref{f.madrasjoinable} are globally Madras joinable.

\begin{lemma}\label{l.strongjoin}
Let $n,m \in 2\N$ and let $\phi^1 \in \sapellward_n$ and $\phi^2 \in \saprward_m$.

Every value 
\begin{equation}\label{e.kvalues}
k \in \Big[ y \big(\righttip(\phi^1)\big)  \, , \,  y\big(\righttip(\phi^1)\big)  +  \min \big\{ n^{1/2}/2,m^{1/2} \big\}  - 1 \Big]
\end{equation}
is such that $\phi^1$ and some horizontal shift of $\phi^2 + ke_2$ is globally Madras joinable. 

Write $\strongjoin_{(\phi^1,\phi^2)}$ for the set of $\vec{u} \in \Z^2$ such that 
the pair $\phi^1$ and $\phi^2 + \vec{u}$ is globally Madras joinable.
Then
$$
  \big\vert \strongjoin_{(\phi^1,\phi^2)} \big\vert \geq   \min \big\{n^{1/2}/2,m^{1/2} \big\}  \, .
$$
\end{lemma}

  \begin{figure}
    \begin{center}
      \includegraphics[width=0.75\textwidth]{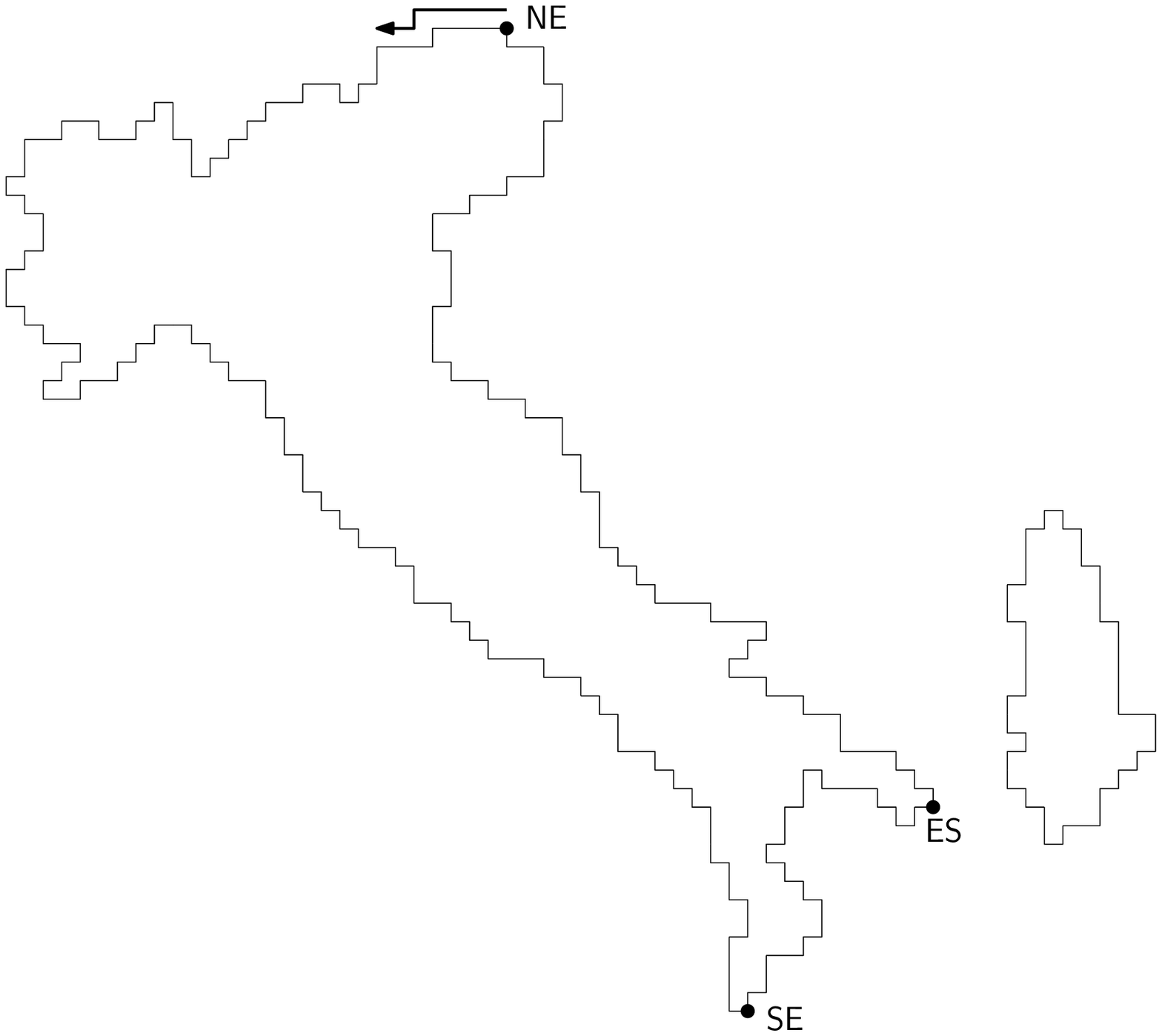}
    \end{center}
    \caption{Illustrating Lemma~\ref{l.strongjoin} with two polygons, $\mathsf{ITA}$ and $\mathsf{ALB}$. Denoting the polygons' lengths by $n$ and $m$, note that the translates of  $\mathsf{ITA}$ and $\mathsf{ALB}$ with $\northeast = 0$ belong to $\sapellward_n$ and $\saprward_m$. 
The start of the counterclockwise tour of $\mathsf{ITA}$ from $\northeast(\mathsf{ITA})$ 
dictated by convention is marked by an arrow.    
    The depicted polygons are globally Madras joinable after a horizontal shift of $\mathsf{ALB}$, with the necessary surgery to $\mathsf{ITA}$ occurring in a vicinity of $\eastsouth(\mathsf{ITA})$.
Vertical shifts of $\mathsf{ALB}$
that leave this polygon in an easterly line of sight from  $\eastsouth(\mathsf{ITA})$ will maintain this state of affairs.
    }
    \label{fig:albita}
  \end{figure}

\noindent{\bf Proof.} Recall that since $\phi^1 \in \sap_n$
and $\phi^2 \in \sap_m$, $\ymax(\phi^1) = \ymax(\phi^2) = 0$.
Note that whenever $k \in \Z$ is such that the two intervals 
$$
\big[ \ymin(\phi^1) , 0 \big] \, \, \, \textrm{and} \, \, \, k + \big[ \ymin(\phi^2) , 0 \big]
$$ 
intersect, there is some horizontal displacement $j \in \Z$ such that $\phi^1$ and $\phi^2 + (j,k)$ are Madras joinable. Note also that $\phi^1 \in \sapellward_n$ satisfies $\ymin(\phi^1)  \leq - n^{1/2}$. Choices of $k \in [ - n^{1/2}  ,  - 1 ]$ thus produce polygon pairs $\big( \phi^1 , \phi^2 + k e_2 \big)$ whose first element contains vertices that are more northerly than any of the second, and which are Madras joinable after a horizontal shift of the second element of the pair. Moreover,  since $\phi^1 \in \sapellward_n$, $y\big( \righttip(\phi^1) \big)$ is at most $- n^{1/2}/2$.  Thus, any value of $k$ in~(\ref{e.kvalues}) lies in $[ - n^{1/2}  ,  - 1 ]$. In view of the Definition~\ref{d.globaljoin} of global join plaquette, in order to verify the lemma's first assertion, it remains to verify that the Madras joinable polygon pair associated by a suitable horizontal shift to any given value of $k$ in~(\ref{e.kvalues}) has the property that its second polygon contains all of the most easterly vertices in either of the two polygons in the pair. Here is the reason that this property holds: 
$\phi^2 \in \saprward_m$ implies that $h(\phi^2) \geq m^{1/2}$, and thus,
for such values of $k$, the $y$-coordinate of $\righttip(\phi^1)$
is shared by a vertex in $\phi^2 + ke_2$; from this, we see that, when this right polygon is shifted horizontally to be Madras joinable with $\phi^1$, it is forced to intersect the half-plane bordered on the left by the vertical line through $\righttip(\phi^1)$ and thus to verify the claimed property. 

The second assertion of the lemma is an immediate consequence of the first.  \qed

\subsection{Polygon joining is almost injective: deriving Proposition~\ref{p.polyjoin.newwork}}

We begin
by reducing this result to the next lemma.

\begin{lemma}\label{l.rgjp.newwork}
For any $\maceta > 0$, there is a constant $\Cpoly = \Cpoly(\maceta) > 0$ such that, for $n  \in 2\N$ satisfying $n + 16 \in \highpolynum_{\Cpcp}$, 
$$
  p_{n + 16}  \geq \frac{1}{\Cpoly \log n} \sum_{j \in 2\N \, \cap \, [2^{i-1},2^i]} p_j p_{n-j} (n-j)^{1/2} \, ,
$$
where $i \in \N$ is chosen so that $n \in 2\N \cap [2^{i},2^{i+1}]$. 
\end{lemma}

\noindent{\bf Proof of Proposition~\ref{p.polyjoin.newwork}.} 
The lemma coincides with the proposition when the two instances of $n + 16$ in its statement are replaced by $n$. To infer the proposition from the lemma (with a relabelling of the value of $\Cpoly(\maceta)$), it is thus enough to establish two things.
First, for some constant $c > 0$, $p_n \geq c \, p_{n+16}$ whenever $n \in 2\N$.  Second, that for $\maceta,\varphi > 0$, there exists $n_0(\maceta,\varphi)$ such that if $n \in 2\N$ satisfies  $n \geq n_0$ and $n + 16 \in \highpolynum_\maceta$, then
$n \in \highpolynum_{\maceta + \varphi}$.
To verify the first statement, recall from \cite[Theorem 7.3.4(c)]{MS93} that $\lim_{n \in 2\N} p_{n+2}/p_n = \mu^2$. Thus, $p_n \geq \mu^{-16} p_{n+16}/2$ for all $n$ sufficiently high; the value of $c > 0$ may be decreased from $\mu^{-16}/2$ if necessary in order that all $n \in 2\N$ satisfy the stated bound.
Given that the first statement holds, 
the second is confirmed by noting that if $n+16 \in \highpolynum_{\maceta}$, then $p_{n+16} \geq (2n)^{-\maceta} \mu^{n+16}$ whenever $n \in 2\N$ satisfies $n \geq 16$.  \qed

\medskip

\noindent{\bf Proof of Lemma~\ref{l.rgjp.newwork}.}
Set 
$A = \bigcup \, \sapellward_j \times \saprward_{n-j}$
where the union is taken over indices 
$j \in 2\N \, \cap \, [2^{i-1},2^i]$ for which $n - j$ is positive. Set $B = \sap_{n + 16}$.

We construct a multi-valued map $\Psi:A \to \mathcal{P}(B)$.
Consider a generic domain point 
 $(\phi^1,\phi^2) \in \sapellward_j \times \saprward_{n-j}$ where $j \in 2\N \, \cap \, [2^{i-1},2^i]$ satisfies $n - j > 0$. 
This point's image under $\Psi$ is defined to be 
 the collection of length-$(n+16)$ polygons 
formed by Madras joining $\phi^1$ and $\phi^2 + \vec{u}$ as $\vec{u}$ ranges over the set $\strongjoin_{(\phi^1,\phi^2)}$ specified in Lemma~\ref{l.strongjoin}.

 Since $j \geq 2^{i-1}$ and $n \leq 2^{i+1}$, we have that $j \geq (n-j)/3$.
By Lemma~\ref{l.strongjoin},
$$
  \Big\vert  \Psi\Big( (\phi^1,\phi^2) \Big) \Big\vert   =  \big\vert \strongjoin_{(\phi^1,\phi^2)} \big\vert \geq \min \big\{ j^{1/2}/2 , (n-j)^{1/2}  \big\} \geq \tfrac{1}{2} \, 3^{-1/2} (n-j)^{1/2} \, .
$$

Applying Lemma~\ref{l.polysetbound.newwork}, we learn that the number of arrows in $\Psi$ is at least
\begin{equation}\label{e.totalarrow}
3^{-1/2}  \, 2^{-4} \sum_{j = 2^{i-1}}^{2^i} p_j p_{n-j}  (n-j)^{1/2} \, ,
\end{equation}
where note that the summand is zero if $j$ is odd. 
 

For a constant $\Ch > 0$, denote by 
$$
H_{n + 16} = \Big\{ \phi \in \sap_{n + 16} : \big\vert \globaljoin_\phi \big\vert \geq  \Ch \log (n+16) \Big\} 
$$
the set of length-$(n + 16)$ polygons with a {\em high} number of global join plaquettes.

We now fix the value of $\Ch$ with which we will work.
Recalling that $n + 16 \in \highpolynum_\maceta$, we may apply Corollary~\ref{c.globaljoin} with $\Cgjone = 1$ 
and the role of $n$ played by $n + 16$.
Set $\Ch > 0$ equal to the value of $\Cgjtwo$ determined via the corollary by this choice of $\Cgjone$ and the value of~$\maceta$.  Since $\psap_{n+16} \big( \vert \globaljoin_\Gamma \vert \geq \Ch \log (n+16) \big) =  \vert H_{n+16} \vert  p_{n + 16}^{-1}$, we find that
\begin{equation}\label{e.hnc}
 \big\vert H_{n + 16} \big\vert \leq \Cgjthree (n + 16)^{-1}  p_{n + 16}  \, .
\end{equation}

The set of preimages under $\Psi$ of a given 
$\phi \in \sap_{n+16}$ may be indexed by the junction plaquette 
associated to the Madras join polygon $J(\phi^1,\phi^2)$ that equals $\phi$. Recalling Definition~\ref{d.globallymj}, this plaquette is a global join plaquette of $\phi$. Thus,  
\begin{equation}\label{e.phisap}
\vert \Psi^{-1}(\phi) \vert \leq \vert \globaljoin_\phi \vert
 \, \, \,
\textrm{for any $\phi \in \sap_{n + 16}$} \,  
.
\end{equation}
 Since $\vert \globaljoin_\phi \vert \leq n + 16$ for all $\phi \in \sap_{n + 16}$, we find that
\begin{equation}\label{e.maxsap}
 \max \big\{ \big\vert \Psi^{-1}(\phi) \big\vert : \phi \in \sap_{n + 16} \big\} \leq n + 16 \, .
\end{equation}

The proof of Lemma~\ref{l.rgjp.newwork} will be completed by considering two cases.

\medskip

 In the first case,
\medskip
\begin{itemize}
\item at least one-half of the arrows in $\Psi$ point to elements of $H_{n+16}$ \, ;
\end{itemize}
\medskip
in the second case, then, at least one-half of these arrows point to elements of
$\sap_{n + 16} \setminus H_{n + 16}$.

\medskip

\noindent{{\bf The first case.}} The inequality
\begin{equation}\label{e.firstcase}
 \big\vert H_{n + 16} \big\vert \, \cdot \, \max \Big\{ \vert \Psi^{-1}(\phi) \vert: \phi \in H_{n + 16} \Big\} \geq \tfrac{1}{2} \cdot  3^{-1/2}  \, 2^{-4} \sum_{j = 2^{i-1}}^{2^i} p_j \, p_{n-j}  (n-j)^{1/2} 
\end{equation}
holds because the left-hand side is an upper bound on the number of arrows in $\Psi$ arriving in $H_{n + 16}$, which is at least one-half of the total number of arrows; and the latter quantity is at least the right-hand side. Applying~(\ref{e.hnc}) and~(\ref{e.maxsap}),
$$
 p_{n + 16} \geq  3^{-1/2} \, 2^{-5} \Cgjthree^{-1}     \sum_{j = 2^{i-1}}^{2^i}  p_j p_{n-j}(n-j)^{1/2}   \, ,
$$ 
so that Lemma~\ref{l.rgjp.newwork} is obtained in the first case.

\medskip

\noindent{{\bf The second case.}}
The quantity
$$
 \Big\vert \, \sap_{n + 16}  \setminus H_{n+16} \, \Big\vert \, \cdot \, \max \Big\{  \big\vert \Psi^{-1}(\phi) \big\vert : \phi \in  \sap_{n + 16}  \setminus H_{n+16} \Big\} $$
is an upper bound on the number of arrows incoming to $\sap_{n + 16}  \setminus H_{n+16}$. In the case that we now consider, the displayed quantity  is thus in view of the total arrow number lower bound~(\ref{e.totalarrow})
 at least $3^{-1/2}  \, 2^{-5} \sum_{j = 2^{i-1}}^{2^i} p_j p_{n-j}  (n-j)^{1/2}$.

 By~(\ref{e.phisap}), for $\phi \in \sap_{n + 16} \setminus H_{n+16}$, $\vert \Psi^{-1}(\phi) \vert \leq \Ch \log (n+16)$.
 Thus,
\begin{equation}\label{e.casetwo}
 p_{n + 16} \geq \Ch^{-1} \big( \log (n+16) \big)^{-1}  3^{-1/2}  \, 2^{-5} \sum_{j = 2^{i-1}}^{2^i} p_j p_{n-j}  (n-j)^{1/2}
 \, ,
\end{equation}
so that Lemma~\ref{l.rgjp.newwork} is proved in the second case also. \qed

\subsection{Inferring Theorem~\ref{t.polydev} from Proposition~\ref{p.polyjoin.newwork}}

In this subsection, we prove Theorem~\ref{t.polydev}.
The key step is the following, a consequence of Proposition~\ref{p.polyjoin.newwork}.

\begin{proposition}\label{p.hpn.propagate}
For each $\varphi > 0$, there exists $i_0(\varphi) > 0$ such that, for any $\delta > 0$,  $a \in (0,1)$ and $i \in \N$ for which $i \geq i_0(\varphi) + 2 (\log 2)^{-1} \big( 3\varphi^{-1} + 1 \big)  \log a^{-1}$,
the condition that
\begin{equation}\label{e.hyphigh}
 \Big\vert \, \highpolynum_{3/2 - \delta} \cap \big[2^{i-1},2^i \big] \, \Big\vert \, \geq \, a \cdot 2^{i-2}
\end{equation}
implies that
$$
 \Big\vert \, \highpolynum_{3/2 - 2\delta + \varphi} \cap \big[2^i,2^{i+1} \big] \, \Big\vert \, \geq \,  \tfrac{1}{8} a^2 \cdot 2^{i-1} \, .
$$
\end{proposition}
The proposition stipulates an unstable state of affairs should Theorem~\ref{t.polydev} fail, as we now outline.
Let $\delta > 0$ and assume that infinitely many dyadic scales are occupied by $\highpolynum_{3/2 - \delta}$ at a small but uniform fraction $a$. The proposition implies that the next scale up from any one of these scales (excepting finitely many) is occupied at a fraction at least $a^2/8$ by the more stringently specified set $\highpolynum_{3/2 - 3\delta/2}$.
We may iterate this inference over several scales and find that, a few dyadic scales higher from where we began, there is a tiny but positive fraction of indices that actually belong to  $\highpolynum_{\maceta}$ for a negative value of $\maceta$. This index set is of course known to be empty. Thus, the  $\highpolynum_{3/2 - \delta}$-occupancy assumption is found to be invalid, and Theorem~\ref{t.polydev} is proved. 

First, we prove Proposition~\ref{p.hpn.propagate} and then, in Proposition~\ref{p.hpn.contradict} and its proof, we implement a rigorous version of the argument just sketched.

\medskip 

\noindent{\bf Proof of Proposition~\ref{p.hpn.propagate}.}
We will prove the result with the choice 
\begin{equation}\label{e.varphiformula}
  i_0(\varphi) = \max \Big\{ 6, \varphi^{-1} \big( 18 + \tfrac{2}{\log 2} \log \Cpoly \big) , \tfrac{\log(2\pi)}{2\log 2} + \tfrac{1}{\log 2} \big( 4\varphi^{-1} + 3/2 \big) \log \big( 4 \varphi^{-1} + 1 \big)  \Big\} \, ,
\end{equation}
where the constant $\Cpoly$ is determined by  Proposition~\ref{p.polyjoin.newwork} with $\maceta = 3$.

Consider the map ${\rm sum}$ from the product of two copies of
$2\N \cap [2^{i-1},2^i]$
to $2\N \cap [2^i,2^{i+1}]$
that sends $(j,k)$ to $j+k$. 
Borrowing the language that we use for multi-valued maps, we think of the function as being specified by a collection of arrows $(j,k) \to j+k$ with target $j+k$. 
Fixing $\delta > 0$, we call an arrow $(j,k) \to j+k$ {\em high}
if both $j$ and $k$ are elements of 
$\highpolynum_{3/2 - \delta}$. 
We also identify the set $\targetmany \subseteq 2\N \cap [2^i,2^{i+1}]$
of {\rm sum} image points {\em targeted by many high arrows},
$$
 \targetmany = \Big\{ n \in  2\N \cap [2^i,2^{i+1}]: n \, \, \textrm{is the target of at least $\tfrac{1}{8} a^2 2^{i-2}$ high arrows}    \Big\} \, .
$$
We now argue that $\targetmany$ occupies a proportion of at least $a^2/8$ of the {\rm sum} image points.
\begin{lemma}\label{l.targetmany}
Suppose that the hypothesis~(\ref{e.hyphigh})
 holds for a given $a \in (0,1)$.
Then
$$
 \big\vert \targetmany \big\vert \geq \tfrac{1}{8} a^2 \cdot 2^{i-1} \, .
$$ 
\end{lemma}
\noindent{\bf Proof.} 
The  {\rm sum} image set  $2\N \cap [2^i,2^{i+1}]$
has $2^{i-1} + 1 \leq 2^i$ elements. Thus, the number of high arrows targeting elements of  $\targetmany^c$
is at most $2^i \cdot \tfrac{1}{8} a^2 2^{i-2} = a^2 \cdot 2^{2i - 5}$. 

By~(\ref{e.hyphigh}), the total number of high arrows is at least $a^2 \cdot 2^{2i- 4}$.
We learn then that at least half of these arrows target elements of $\targetmany$. However, no image point is the target of more than $2^{i-2} + 1 \leq 2^{i-1}$ arrows. Thus, there must be at least $a^2 \cdot 2^{2i - 5} \cdot 2^{1-i}$ elements of $\targetmany$. This completes the proof. \qed
 
\medskip

The next lemma and its proof concern a useful additive stability property of the high polygon number sets: for $\maceta > 0$, $\highpolynum_{\maceta} \times \highpolynum_{\maceta}$
is mapped by {\rm sum} into  $\highpolynum_{2\maceta}$. This fact is useful when we seek to apply Proposition~\ref{p.polyjoin.newwork}: the proposition will be useful when the right-hand side in~(\ref{e.polyjoin.newwork}) is in a suitable sense large, but, in this circumstance, the $\highpolynum$-stability property  verifies the proposition hypothesis that $n$ is an element of some $\highpolynum$-set.
 
\begin{lemma}\label{l.target.hpn}
If $n \in 2\N \cap [2^i,2^{i+1}]$ is the target of a high arrow, then~$n \in \highpolynum_{3 - 2\delta}$.
\end{lemma}
\noindent{\bf Proof.}
Let $(j,n-j)$ be the high arrow that targets $n$. 
By~(\ref{e.threepbound}), $p_n \geq p_j p_{n-j}$.
Since $j$ and $n-j$ are elements of $\highpolynum_{3/2 - \delta}$, we have that 
\begin{equation}\label{e.mintwop}
\min \{ p_j \mu^{-j}, p_{n-j} \mu^{-(n-j)} \} \geq n^{-3/2 + \delta}  \, ,
\end{equation}
where we also used $\delta \leq 3/2$, which bound follows from $\highpolynum_{3/2 - \delta} \not= \emptyset$ and $p_n \leq \mu^n$ from~(\ref{e.pbasicthree}).
 Thus, $p_n \geq n^{-3 + 2\delta} \mu^{2n}$. \qed

\medskip

Fix an element $n \in \targetmany$. Lemma~\ref{l.target.hpn}
shows that Proposition~\ref{p.polyjoin.newwork} is applicable with $\maceta = 3$. Thus,
\begin{equation}\label{e.pnlb}
 p_n \geq \frac{1}{\Cpoly \log n} \sum (n-j)^{1/2} p_j \, p_{n-j}  \, ,
\end{equation}
where the sum is taken over $j \in 2\N \, \cap \, [2^{i-1},2^i]$
such that $(j,n-j) \to n$ is a high arrow.

By its definition and the structure of the map ${\rm sum}$,  the set $\targetmany$ is in fact a subset of the even integers in the interval $[2^i + A, 2^{i+1} - A]$, where $A = \tfrac{1}{4} a^2 2^{i-2} - 2$.

Recalling (\ref{e.varphiformula}),
note that the hypotheses of the proposition entail that 
$i \geq 6 + 2(\log 2)^{-1} \log a^{-1}$.
Thus, $a \geq 2^{3 - i/2}$.
Let $j \in [2^{i-1},2^i]$.   
Since $n \in \targetmany$, we see from $n - j \geq A$
and this lower bound on $a$ that $n-j$ is at least
$a^2 2^{i-5}$.

Note then that
\begin{eqnarray*}
 p_n & \geq &  \frac{1}{\Cpoly \log n} a 2^{(i-5)/2} \sum  p_j \, p_{n-j} \\
  & \geq &     \frac{1}{\Cpoly \log n} a 2^{(i-5)/2}   \cdot  \tfrac{1}{8} a^2 2^{i-2}  \cdot n^{-3 + 2 \delta} \mu^{n} \\
  & \geq & 
   \frac{1}{2^{5/2+3+2} \Cpoly \log n} a^3 \cdot 2^{-3/2} n^{3/2} \cdot n^{-3 + 2 \delta} \mu^{n} = \frac{1}{2^{9} \Cpoly \log n} a^3  \cdot n^{-3/2 + 2 \delta} \mu^{n} \, .  
\end{eqnarray*}
In the first line, we sum over the same set of values of $j$ as we did in~(\ref{e.pnlb}). The inequality in this line is due to
(\ref{e.pnlb}) and $n - j \geq a^2 2^{i-5}$.
The second inequality is due to $n \in \targetmany$ and the bound~(\ref{e.mintwop}).  The final line arises because  $2^{i} \geq n/2$.

Finally, we argue that 
\begin{equation}\label{e.varphiargument}
 2^{-9} \Cpoly^{-1} a^3 \big(\log n \big)^{-1} \geq n^{-\varphi} \, .
\end{equation}
Combining with the last display, we learn that $p_n \geq n^{-3/2 + 2\delta - \varphi} \mu^n$, which is to say, $\thet_n \leq 3/2 - 2\delta + \varphi$. Recalling that this inference has been made for any element $n \in \targetmany$, we note that Lemma~\ref{l.targetmany}
completes the proof of Proposition~\ref{p.hpn.propagate}. 

It remains to confirm~(\ref{e.varphiargument}). This follows from two estimates: first, $\log n \leq n^{\varphi/2}$;
second, $2^{-9} \Cpoly^{-1} a^3 \geq n^{-\varphi/2}$.
The first is these follows by taking $k = \lceil 4\varphi^{-1} \rceil$ in $e^{n^{\varphi/2}} \geq \tfrac{1}{k!} n^{k \varphi/2}$
using $k! \leq (2\pi)^{1/2} \big( 4 \varphi^{-1} + 1 \big)^{4 \varphi^{-1} + 3/2}$ as well as $n \geq 2^i$ and $i \geq i_0(\varphi)$. The second follows from $n \geq 2^i$
and the hypothesised lower bound on~$i$. \qed

\begin{proposition}\label{p.hpn.contradict}
Let $\delta > 0$ and  $a \in (0,1)$. Let $i_0(\delta/2) > 0$ be  specified by Proposition~\ref{p.hpn.propagate}. 
Suppose that $i \in \N$ satisfies 
$$
i \geq 
 i_0(\delta/2) +  f(\delta,a) + 2 (\log 2)^{-1} \big( 6 \delta^{-1} + 1 \big) \log a^{-1} \, , 
$$ 
where  
$$
f(\delta,a) = 4 \big( \log 2 \big)^{-1}  \big( 6 \delta^{-1} + 1 \big)  \delta^{-( \log 3/2 )^{-1} \log 2}   \log \big( 8 a^{-1} \big) \, .
$$
Then
$$
 \Big\vert \, \highpolynum_{3/2 - \delta} \cap \big[2^{i-1},2^i \big] \, \Big\vert \, < \, a \cdot 2^{i-2} \, .
$$
\end{proposition}
\noindent{\bf Proof of Theorem~\ref{t.polydev}.} This is an immediate consequence of the proposition. \qed

\medskip

\noindent{\bf Proof of Proposition \ref{p.hpn.contradict}.}
Suppose the contrary. Let $\delta > 0$, $a \in (0,1)$ and 
\begin{equation}\label{e.imin}
i \geq i_0 \big( \delta/2 \big) \, + \, f(\delta,a) \, + \, 2 (\log 2)^{-1} \big( 6 \delta^{-1} + 1 \big)  \log a^{-1} 
\end{equation}
be such that~(\ref{e.hyphigh}) holds.

For $k \in \N$, set $h_k(a) = 8^{1 - 2^k} a^{2^k}$, and note for future reference that the recursion $h_{k+1}(a) = h_k(a)^2/8$ is satisfied with initial condition $h_0(a) = a$.

Note that the inequality 
\begin{equation}\label{e.hkdelta}
  2 (\log 2)^{-1} \big( 6 \delta^{-1} + 1 \big) \log h_k(a)^{-1} > f(\delta,a) + 2 (\log 2)^{-1}  \big( 6 \delta^{-1} + 1 \big)  \log a^{-1}  
\end{equation} 
  is satisfied for all sufficiently high $k \in \N$. Let $K \in \N$
  denote the smallest value of $k \in \N$ that satisfies this bound. Our choice of the function $f$ entails that $K > 1 + \big( \log 3/2 \big)^{-1}  \log \delta^{-1}$.

We claim that, for any $k \in [0,K]$,
\begin{equation}\label{e.k.hpn}
 \Big\vert \, \highpolynum_{3/2 - (3/2)^k \delta} \cap \big[2^{i+k-1},2^{i+k} \big] \, \Big\vert \, \geq \,  h_k(a) \cdot 2^{i-k-2} \, .
\end{equation}
This bound may be proved by induction on $k$, where the
base case $k=0$ is validated by our assumption that~(\ref{e.hyphigh}) holds.
Proposition~\ref{p.hpn.propagate} yields the bound for $k \in [1,K]$ as we now explain.
Suppose that the bound holds for index $k-1$. We seek to apply Proposition~\ref{p.hpn.propagate}
with the roles of $\delta$, $\varphi$, $i$ and $a$ being played by $(3/2)^{k-1} \delta$, $(3/2)^{k-1} \delta/2$,  $i + k -1$
and $h_{k-1}(a)$. The conclusion of the proposition 
 is then indeed the bound at index $k$, but we must check that the proposition's hypotheses are valid. 
Beyond the inductively supposed hypothesis, we must confirm that
\begin{equation}\label{e.ikhyp}
 i + k -1 \geq  i_0 ( \varphi ) \, + \, 
  2 (\log 2)^{-1} \big( 3 \varphi^{-1} + 1 \big) \log h_{k-1}(a)^{-1}   \, ,
\end{equation}
with $\varphi =  (3/2)^{k-1} \delta/2$. Note that since $k \in [0,K]$,
the inequality~(\ref{e.hkdelta}) is violated at index $k-1$.
Note that $i + k - 1 \geq i$ may be bounded below by~(\ref{e.imin}) and then, by this violation, by
$i_0(\delta/2) +   2 (\log 2)^{-1} \big( 6 \delta^{-1} + 1 \big) \log h_{k-1}(a)^{-1}$. Since $\varphi \geq \delta/2$ and $i_0$ is a decreasing function of $\varphi > 0$, we do indeed verify~(\ref{e.ikhyp}) and so complete the derivation of~(\ref{e.k.hpn}) for each $k \in [0,K]$. 

Consider now (\ref{e.k.hpn}) with $k = K$:
we see that 
$\big\vert \highpolynum_{\zeta} \cap \big[2^{i+
 K-1},2^{i+K} \big]  \big\vert  \geq   h_K(a) \cdot 2^{i-K-2} > 0$ where $\zeta = 3/2 - (3/2)^K \delta$.
However, by recalling that $K$ exceeds $\big( \log 3/2 \big)^{-1}  \log \delta^{-1}$ by more than one, we find $\zeta$ 
 to be negative, so that $p_n \leq \mu^n$ implies that $\highpolynum_\zeta = \emptyset$.
  By obtaining this contradiction, we have completed the proof of Proposition~\ref{p.hpn.contradict}. \qed

\section{General dimension: deriving Theorem~\ref{t.threed}}\label{s.threed}

Here we adapt the proof of Theorem~\ref{t.polydev} from Section~\ref{s.three}
to prove Theorem~\ref{t.threed} concerning the 
 $3$-edge self-avoiding walk model. 
 
First we recall that in the introduction we asserted the existence of $\muhat =  \lim_{n \in 2\N} \barp_n^{1/n}$. 
In fact, we also have that $ \lim_{n \in \N} \carp_n^{1/n}$
exists and equals this value.  Moreover, we asserted that
$c_n \geq \muhat^n$ for $n \in \N$. 
These claims are counterparts of (1.2.10) and (3.2.9) in \cite{MS93} for the $3$-edge self-avoiding walk model.  The proof of these results, in particular the Hammersley-Welsh unfolding argument needed to show~\cite[(3.2.9)]{MS93}, can be adapted with only minor changes, a discussion of which we omit.

We respecify the index set from Definition~\ref{d.highpolynum}
so that,
for $\maceta > 0$, 
$\highpolynum_\maceta = \big\{ n \in 2\N : \barp_n \geq n^{- \maceta} \muhat^n \big\}$.
With this change made, we may state an analogue of Proposition~\ref{p.polyjoin.newwork}.
\begin{proposition}\label{p.polyjoin3d}
Let $d \geq 2$. For any $\maceta > 0$, there is a constant $\Cpoly = \Cpoly(\maceta) > 0$ such that, for $n \in 2\N \cap \highpolynum_{\Cpcp}$,
$$
 \barp_n \geq \frac{1}{\Cpoly \log n} \sum_{j \in 2\N \, \cap \, [2^{i-1},2^i]} (n-j)^{1 - 1/d} \barp_j \, \barp_{n-j}  \, ,
$$
where $i \in \N$ is chosen so that $n \in 2\N \cap [2^i,2^{i+1}]$. 
\end{proposition}

Our proof follows closely the structure of Theorem~\ref{t.polydev}'s, with this section having five subsections each of which plays a corresponding role  to its counterpart in Section~\ref{s.three}. As such, the first four subsections lead to the proof of Proposition~\ref{p.polyjoin3d}, and, in the final one, Theorem~\ref{t.threed} is proved from the proposition.

\subsection{Specifying the joining of two polygons}

\begin{definition}
Let $\phi$ be a polygon.  A {\em join} edge of $\phi$ is a nearest neighbour edge in $\Z^d$ that is traversed precisely two times by $\phi$, with one crossing in each direction. 
\end{definition}
Note that when a join edge is removed from a polygon, two polygons result.

\begin{definition}
For $I$ any $d-1$ element subset of $\{1,2,\cdots,d\}$, let $\proj_I:\Z^d \to \Z^d$ denote the projection onto the axial plane containing the vectors $e_i$ for $i \in I$. Write $\proj = \proj_{2,\cdots,d}:\Z^d \to \{ 0 \} \times \Z^{d-1}$. 
\end{definition}
We present an analogue of Madras' local surgery procedure for polygon joining, respecifying the join $J(\phi,\phi')$ of two polygons. 
For $n,m \in 2\N$, let $\phi$ and $\phi'$ be polygons of lengths $n$ and $m$. Suppose that $\proj(\phi) \cap \proj(\phi') \not= \emptyset$. 
If $\phi'$ is translated to a sufficiently high coordinate in the $e_1$-direction,   the vertex sets of $\phi$ and the $\phi'$-translate will be disjoint. Make such a translation and then translate $\phi'$ backwards in the negative $e_1$-direction  until the final location at which the translate's vertices remain disjoint from those of $\phi$. When two polygons have such a relative position, they are in a situation comparable to being Madras joinable; we call them {\em simply} joinable. Note that there exists an $e_1$-oriented nearest neighbour edge of~$\Z^d$
with one endpoint in one of the polygon vertex sets and the other in the other such set.
Let $e$ denote the maximal edge among these (according to some fixed ordering of nearest neighbour edges of~$\Z^d$). Define $J(\phi,\phi')$ to be the length~$n + m + 2$ polygon formed by following the trajectory of $\phi$ until an endpoint of $e$ is encountered, crossing $e$, following the whole trajectory of $\phi'$, recrossing $e$, and completing the trajectory of~$\phi$. The edge $e$ will be called the {\em junction} edge in this construction.

Let the left vertex $\leftv(\phi)$ of a polygon $\phi$ be the lexicographically minimal vertex in $\phi$ (so that $\leftv(\phi)$ has minimal $e_1$-coordinate). For our purpose, the left vertex will play a counterpart role  to that of the northeast vertex  in the setting of the proof of Theorem~\ref{t.polydev}. For $n \in 2\N$, we define $\saptd_n$ to be the set of polygons of length~$n$ whose left vertex is the origin. Note that $\barp_n = \# \saptd_n$. We adopt a convention for parametrizing any $\phi \in \saptd_n$, so that we may write $\phi:[0,n] \to \Z^d$.
The polygon $\phi$ visits the origin, possibly more than once.
For each such visit, we may record two lists of length $n$ of the vertices consecutively visited by $\phi$, with one having the opposite orientation of time to the other. 
Among these finitely many lists, we select the one of minimal lexicographical order, and take $\phi:[0,n] \to \Z^d$ to be this list.

We also define the right vertex $\rightv(\phi)$ of a polygon $\phi$ be the  lexicographically maximal vertex in $\phi$ (which is thus one of maximal $e_1$-coordinate).

\subsection{Global join edges are sparse}

Here we present analogues to Definition~\ref{d.globaljoin}, Proposition~\ref{p.globaljoin} and this result's proof.

\begin{definition}
Let $n \in 2\N$ and let $\phi \in \saptd_n$.  A join edge $e$ of $\phi$ is called {\em global} if the two polygons comprising $\phi$ without $e$ may be labelled  $\phi^\ell$ and $\phi^r$ in such a way that 
\begin{itemize}
\item every vertex of maximal $e_1$-coordinate in $\phi^\ell \cup \phi^r$ belongs to $\phi^r$;
\item and every vertex of minimal $e_1$-coordinate in $\phi^\ell \cup \phi^r$ belongs to $\phi^\ell$.
\end{itemize}
Write $\globaljoin_\phi$ for the set of global join edges of the 
polygon~$\phi$.
\end{definition}

Proposition~\ref{p.globaljoin}'s most direct analogue holds.
\begin{proposition}\label{p.globaljoin3d}
There exists $c > 0$ such that, for $n \in 2\N$ and any $k \in \N$,
$$
  \# \Big\{ \phi \in \saptd_{n} : \big\vert \globaljoin_\phi \big\vert \geq k \Big\} \leq  c^{-1} 2^{-k \muhat^{-2}/2}    \muhat^n   \, .
$$
\end{proposition}
\noindent{\bf Proof.} Let $\phi \in \saptd_n$. The right vertex $\rightv(\phi)$ replaces $\eastsouth(\phi)$ in this proof. 
Set $j \in [0,n]$ so that $\phi_j = \rightv(\phi)$; write $\phi^1 =  \phi_{[0,j]}$ and $\phi^2 = \phi_{[j,n]}$, and denote by $\reflect_z$ reflection in the hyperplane with normal vector~$e_1$ that passes through $z \in \Z^d$. Then set
$$
\mathscr{S}(\phi) = \phi^1 \circ \reflect_{\rightv(\phi)}(\phi^2) \, .
$$
The proof proceeds as Proposition \ref{p.globaljoin}'s did.
In the present case, $\mathscr{S}$ is mapping $\saptd_n$
into the set of bridges of length~$n$, namely into the set of $3$-edge length-$n$ self-avoiding walks whose minimal and maximal $e_1$-coordinate is attained at the start and the endpoint. In fact, bridges are usually defined so that this attainment is unique at one of the extremes. The addition of an extra edge at the end can assure this uniqueness. Using the classical upper bound on bridge number discussed in the proof of~(\ref{e.rhswalkbound}) (which is easily adapted to the $3$-edge case), we see that the image set $\mathscr{S}(\saptd_n)$
has cardinality at most $\muhat^{n+1}$. 
 
Another modification is made in specifying the alternative walks
$\mathscr{S}_\kappa(\phi)$, where $\kappa \subseteq \{ 1,\cdots,r\}$. Denote the $i\textsuperscript{th}$ global join edge by $g_i \in \globaljoin_\phi$. Because $g_i$ is a join edge, it is traversed twice by~$\phi$. The removal of $g_i$ from $\phi$
results in two polygons $\phi^\ell$ and $\phi^r$. Since $g_i$ is a global join edge, we have that 
$\leftv(\phi) \in \phi^\ell$ and $\rightv(\phi) \in \phi^r$.
Thus, $\phi_{[0,j]}$ contains one traversal of the edge $g_i$ made by $\phi$,
and $\phi_{[j,n]}$ the other. 
In light of this, we may describe the local modication for index~$i$, analogous to the three step subpath around $P^i$ in the original proof. 
Consider the edge in  $\mathscr{S}(\phi)$ that is the reflected image of $g_i$.
The modified walk is specified by insisting that it follows the course of  $\mathscr{S}(\phi)$ until this particular edge is traversed; after the traversal, the new walk performs a two-step move, crossing back over the edge, and then back again; after this, it pursues the remaining trajectory of   $\mathscr{S}(\phi)$.
  Note that no edge is traversed more than three times in the resulting definition of $\mathscr{S}_\kappa(\phi)$, and also that, when the post-concatenation part of this walk is reflected back, the outcome is a walk with an edge being traversed four times associated to any surgery, rendering the location of these surgeries detectable. (This property serves to explain our use of $3$-edge self-avoiding walks.) \qed  
 
\subsection{Left-right polygon pairs}

We now specify a counterpart  $\lrpp_{(k,\ell)}$  to the set $\sapellward_k \times \saprward_\ell$ of pairs of left and right polygons. The definition-lemma-definition-lemma structure of the counterpart Subsection~\ref{s.refpert.newwork} is maintained.


For $i \in [1,d]$, the $e_i$-span of a polygon is defined to be the difference between the maximum and the minimum $e_i$-coordinate of elements of the vertex set of $\phi$. Write  $\xspan(\phi)$ for the $e_1$-span of $\phi$.
 
\begin{definition} 
Let $k,\ell \in 2\N$.
A $(k,\ell)$ left-right polygon pair $(\phi,\phi')$
is an element of the union of $\saptd_k \times \saptd_\ell$
and $\saptd_\ell \times \saptd_k$ such that
\begin{itemize}
\item $\xspan(\phi) \geq \xspan(\phi')$, 
\item and $\vert \proj(\phi') \vert \geq (3d)^{-(1 - 1/d)} l(\phi')^{1 - 1/d}$, where $l(\phi')$ denotes the length of~$\phi'$.
\end{itemize}
Write $\lrpp_{(k,\ell)}$ for the set of  $(k,\ell)$ left-right polygon pairs.
\end{definition}

\begin{lemma}\label{l.polysetbound3d}
For $k,\ell \in 2\N$,
$$
\Big\vert \, \lrpp_{(k,\ell)} \,  \Big\vert \, \geq \, \tfrac{1}{2} d^{-2} \cdot \big\vert \saptd_k \big\vert \cdot  \big\vert \saptd_\ell \big\vert   
\, .
$$
\end{lemma}
\noindent{\bf Proof.}
Begin with a polygon pair $(\phi,\phi') \in \saptd_k \times \saptd_\ell$. By relabelling and reordering the pair, we may ensure that $\xspan(\phi)$
is at least the span of $\phi'$ in any of the~$d$ coordinates.
This relabelling is responsible for a factor of $(2d)^{-1}$ on the right-hand side. Now further relabel the coordinate axes in regard only to $\phi'$ in order that the cardinality of $\proj(\phi')$ for the relabelled $\phi'$
is at least as large as the cardinality associated to the other $d-1$ directions. This entails the appearance of a further factor of $d^{-1}$ on the right-hand side. It is easily seen that $\vert V(\phi') \vert \geq \tfrac{1}{3d}l(\phi')$ (since $\phi'$ is a $3$-edge self-avoiding polygon), and, as such, the lemma will be proved if we can show that the maximum cardinality among the axial projections of $\proj(\phi')$ is at least $\vert V(\phi') \vert^{1 - 1/d}$.
 This is a consequence of the Loomis-Whitney inequality~\cite{LoomisWhitney}: for any $d \geq 2$ and any finite $A \subseteq \Z^d$, the product of the cardinality of the projections of $A$ onto each of the $d$ axial hyperplanes is at least $\vert A \vert^{d-1}$. Also using the arithmetic-geometric mean inequality, one of the projections has cardinality at least $\vert A \vert^{1 - 1/d}$. \qed

\medskip

For polygons $\phi$ and $\phi'$, write $\strongtdjoin_{(\phi,\phi')}$ for the set of $\vec{u} \in \Z^d$ such that the pair 
$\big( \phi,\phi' + \vec{u} \big)$
is {\em globally joinable}, which is to say that 
\begin{itemize}
\item
$\big( \phi,\phi' + \vec{u} \big)$ is simply joinable;
\item any element of minimal $e_1$-coordinate  among the vertices of $\phi$ or $\phi' + \vec{u}$ is a vertex in $\phi$;
\item and any element of maximal $e_1$-coordinate  among the vertices of $\phi$ or $\phi' + \vec{u}$ is a vertex in $\phi' + \vec{u}$.
\end{itemize}

\begin{lemma}\label{l.claimmin}
For 
$k,\ell \in 2\N$, if 
$(\phi,\phi') \in \lrpp_{(k,\ell)}$, 
then 
$$
\big\vert \strongtdjoin_{(\phi,\phi')} \big\vert \geq (3d)^{-(1 - 1/d)} \min \big\{ k^{1-1/d} , \ell^{1-1/d} \big\} \, .
$$
\end{lemma}
\noindent{\bf Proof.}
Suppose that  $\vec{u} \in \{ 0 \} \times \Z^{d-1}$ is such that  $\proj(\phi' + \vec{u})$ contains the vertex $\proj\big(\rightv(\phi)\big)$. Due to our choice of $\phi'$, such $\vec{u}$ number 
at least the right-hand side in the inequality in the lemma's statement.  It is thus sufficient to argue that, for any such $\vec{u}$, there exists $k \in \Z$ for which the pair $(\phi,\phi' + \vec{u} + k e_1)$ is globally joinable. Recall that there is a maximal choice of $k \in \Z$ such that the polygon pair   $\big(\phi,\phi' + \vec{u} + (k-1) e_1 \big)$ is not vertex disjoint but the pair  $\big(\phi,\phi' + \vec{u} + k e_1 \big)$ is vertex disjoint, and that the pair  $\big(\phi,\phi' + \vec{u} + k e_1 \big)$ is said to be simply joinable. Specifying $k \in \Z$ in this way, we must also check that  $\vec{u} + k e_1$ verifies the second and third conditions for membership of~$\strongtdjoin_{(\phi,\phi')}$. To do  this, note that 
the polygon $\phi' + \vec{u} + k e_1$ contains a vertex with an $e_1$-coordinate strictly exceeding that of any vertex in $\phi$, so that the third condition is confirmed. Moreover, this same fact implies the second condition, because   the $e_1$-span of $\phi$ exceeds that of $\phi'$. We have proved Lemma~\ref{l.claimmin}. \qed

\subsection{Proof of Proposition~\ref{p.polyjoin3d}}
We begin with Lemma~\ref{l.rgjp.newwork}'s counterpart.
\begin{lemma}\label{l.rgjp.newwork.hd}
For any $\maceta > 0$, there is a constant $\Cpoly = \Cpoly(\maceta) > 0$ such that, for $n  \in 2\N$ satisfying $n + 2 \in \highpolynum_{\Cpcp}$, 
$$
  \hat{p}_{n + 2}  \geq \frac{1}{\Cpoly \log n} \sum_{j \in 2\N \, \cap \, [2^{i-1},2^i]} \hat{p}_j \hat{p}_{n-j} (n-j)^{1 - 1/d} \, ,
$$
where $i \in \N$ is chosen so that $n \in 2\N \cap [2^{i},2^{i+1}]$. 
\end{lemma}
Lemma~\ref{l.rgjp.newwork.hd}  implies Proposition~\ref{p.polyjoin3d} as Lemma~\ref{l.rgjp.newwork} did Proposition~\ref{p.polyjoin.newwork}, with the assertion counterpart to \cite[Theorem 7.3.4(c)]{MS93} that $\lim_{n \in 2\N} \hat{p}_{n+2}/\hat{p}_n = \muhat^2$ being used. 
The cited result has a non-trivial proof, but the changes needed to the proof are trivial. We do not provide details, but mention very briefly the constructs that would be used:
a pattern is a local configuration that is capable of appearing as a fragment in the middle of a long walk; a pair of patterns called type $I$ and $II$ may be specified so that the second is formed from the first by a deformation that adds a net total of two edges;   Kesten's pattern theorem~\cite[Theorem 1]{kestenone} 
concerning the ubiquity in any typical long walk of a given pattern must be rederived for $3$-edge self-avoiding walks, and then applied to argue that
changing one uniformly chosen type~$II$ pattern to a type~$I$ pattern in an element of $\psaptd_{n+2}$ generates a law on $\saptd_n$ which is provably close to the uniform law $\psaptd_n$.

\noindent{\bf Proof of Lemma~\ref{l.rgjp.newwork.hd}.} This is similar to the proof of Lemma~\ref{l.rgjp.newwork}. We consider the multi-valued map $\Psi:A \to \mathcal{P}(B)$, where 
$$
A = \bigcup_{j \in 2\N \, \cap \, [2^{i-1},2^i]} \lrpp_{(n-j,j)} 
$$
and $B = \saptd_{n+2}$. The map $\Psi$
associates to each $(\phi^1,\phi^2) \in  \lrpp_{(n-j,j)}$, $j \in 2\N \, \cap \, [2^{i-1},2^i]$, the set of length-$(n+2)$ polygons 
formed by simply joining $\phi^1$ and $\phi^2 + \vec{u}$ for choices of $\vec{u}$ in $\strongtdjoin_{(\phi^1,\phi^2)}$.
That the image set may be chosen to be $B$ is due to the left vertex of any formed polygon being the origin (as the second property in the definition of globally joinable indicates).

The proof follows the earlier one. In place of the assertion leading to~(\ref{e.phisap}) that the junction plaquette associated to the join polygon of a globally Madras joinable polygon pair is a global join plaquette, we instead assert that  in the join polygon $J(\phi^1,\phi^2)$ of two concerned polygons $\phi^1$ and $\phi^2$, the junction edge is a global join edge (which statement is a trivial consequence of the definition of globally joinable). Note that Corollary~\ref{c.globaljoin} may be invoked with obvious notational changes, because Proposition~\ref{p.globaljoin3d} replaces Proposition~\ref{p.globaljoin}.  \qed

\subsection{Proof of Theorem~\ref{t.threed}}

Since Proposition~\ref{p.polyjoin3d}
varies from 
Proposition~\ref{p.polyjoin.newwork} 
only by the relabelling of some notation,
Theorem~\ref{t.threed} 
follows from the former proposition as Theorem~\ref{t.polydev}
did from the latter after evident notational changes are made.
There are two manifestations of the value of dimension in these arguments, in the proof of Lemma~\ref{l.target.hpn} and at the end of the proof of Proposition~\ref{p.hpn.contradict}, when $p_n \leq (d-1)\mu^n$ is applied for $d=2$. The proofs now work provided that the dyadic scale $i$ is sufficiently high. \qed

\bibliographystyle{plain}

\bibliography{saw}

\begin{thebibliography}{10}

\bibitem{BBSlog}
R.~Bauerschmidt, D.~Brydges, and G.~Slade.
\newblock Logarithmic correction for the susceptibility of the 4-dimensional
  weakly self-avoiding walk: a renormalisation group analysis.
\newblock arXiv:1403.7422. Commun. Math. Phys., to appear, 2014.

\bibitem{BDGS11}
R.~Bauerschmidt, H.~Duminil-Copin, J.~Goodman, and G.~Slade.
\newblock Lectures on self-avoiding walks.
\newblock In D.~Ellwood, C.~Newman, V.~Sidoravicius, and W.~Werner, editors,
  {\em Lecture notes, in Probability and Statistical Physics in Two and More
  Dimensions}. CMI/AMS -- Clay Mathematics Institute Proceedings, 2011.

\bibitem{BSRigRen3}
Roland Bauerschmidt, David~C. Brydges, and Gordon Slade.
\newblock A {R}enormalisation {G}roup {M}ethod. {III}. {P}erturbative
  {A}nalysis.
\newblock {\em J. Stat. Phys.}, 159(3):492--529, 2015.

\bibitem{BMB09}
Mireille Bousquet-M{\'e}lou and Richard Brak.
\newblock Exactly solved models.
\newblock In {\em Polygons, polyominoes and polycubes}, volume 775 of {\em
  Lecture Notes in Phys.}, pages 43--78. Springer, Dordrecht, 2009.

\bibitem{BSRigRen1}
David~C. Brydges and Gordon Slade.
\newblock A {R}enormalisation {G}roup {M}ethod. {I}. {G}aussian {I}ntegration
  and {N}ormed {A}lgebras.
\newblock {\em J. Stat. Phys.}, 159(3):421--460, 2015.

\bibitem{BSRigRen2}
David~C. Brydges and Gordon Slade.
\newblock A {R}enormalisation {G}roup {M}ethod. {II}. {A}pproximation by
  {L}ocal {P}olynomials.
\newblock {\em J. Stat. Phys.}, 159(3):461--491, 2015.

\bibitem{BSRigRen4}
David~C. Brydges and Gordon Slade.
\newblock A {R}enormalisation {G}roup {M}ethod. {IV}. {S}tability {A}nalysis.
\newblock {\em J. Stat. Phys.}, 159(3):530--588, 2015.

\bibitem{BSRigRen5}
David~C. Brydges and Gordon Slade.
\newblock A {R}enormalisation {G}roup {M}ethod. {V}. {A} {S}ingle
  {R}enormalisation {G}roup {S}tep.
\newblock {\em J. Stat. Phys.}, 159(3):589--667, 2015.

\bibitem{ontheprob}
Hugo Duminil-Copin, Ioan Manolescu, Alexander Glazman, and Alan Hammond.
\newblock On the probability that self-avoiding walks ends at a given point.
\newblock arXiv:1305.1257. Ann. Probab., to appear, 2013.

\bibitem{Dup89}
B.~Duplantier.
\newblock Fractals in two dimensions and conformal invariance.
\newblock {\em Phys. D}, 38(1-3):71--87, 1989.
\newblock Fractals in physics (Vence, 1989).

\bibitem{Dup90}
B.~Duplantier.
\newblock Renormalization and conformal invariance for polymers.
\newblock In {\em Fundamental problems in statistical mechanics {VII}
  ({A}ltenberg, 1989)}, pages 171--223. North-Holland, Amsterdam, 1990.

\bibitem{Flory}
P.~Flory.
\newblock {\em Principles of Polymer Chemistry}.
\newblock Cornell University Press, 1953.

\bibitem{CJC}
Alan Hammond.
\newblock On self-avoiding polygons and walks: counting, joining and closing.
\newblock arXiv:1504.05286, 2017.

\bibitem{snakemethodpattern}
Alan Hammond.
\newblock On self-avoiding polygons and walks: the snake method via pattern
  fluctuation.
\newblock \url{math.berkeley.edu/~alanmh/papers/snakemethodpattern.pdf}, 2017.

\bibitem{snakemethodpolygon}
Alan Hammond.
\newblock On self-avoiding polygons and walks: the snake method via polygon
  joining.
\newblock \url{math.berkeley.edu/~alanmh/papers/snakemethodpolygon.pdf}, 2017.

\bibitem{HS92}
T.~Hara and G.~Slade.
\newblock Self-avoiding walk in five or more dimensions. {I}. {T}he critical
  behaviour.
\newblock {\em Comm. Math. Phys.}, 147(1):101--136, 1992.

\bibitem{HS92b}
Takashi Hara and Gordon Slade.
\newblock The lace expansion for self-avoiding walk in five or more dimensions.
\newblock {\em Rev. Math. Phys.}, 4(2):235--327, 1992.

\bibitem{kestenone}
H.~Kesten.
\newblock On the number of self-avoiding walks.
\newblock {\em J. Mathematical Phys.}, 4:960--969, 1963.

\bibitem{Lawler13}
G.~Lawler.
\newblock Random walk problems motivated by statistical physics.
\newblock http://www.math.uchicago.edu/~lawler/russia.pdf, 2013.

\bibitem{LoomisWhitney}
L.~H. Loomis and H.~Whitney.
\newblock An inequality related to the isoperimetric inequality.
\newblock {\em Bull. Amer. Math. Soc}, 55:961--962, 1949.

\bibitem{Madras95}
N.~Madras.
\newblock A rigorous bound on the critical exponent for the number of lattice
  trees, animals, and polygons.
\newblock {\em Journal of Statistical Physics}, 78(3-4):681--699, 1995.

\bibitem{MS93}
N.~Madras and G.~Slade.
\newblock {\em The self-avoiding walk}.
\newblock Probability and its Applications. Birkh\"auser Boston Inc., Boston,
  MA, 1993.

\bibitem{Madras91}
Neal Madras.
\newblock Bounds on the critical exponent of self-avoiding polygons.
\newblock In {\em Random walks, {B}rownian motion, and interacting particle
  systems}, volume~28 of {\em Progr. Probab.}, pages 359--371. Birkh\"auser
  Boston, Boston, MA, 1991.

\bibitem{Nie82}
B.~Nienhuis.
\newblock Exact critical point and critical exponents of ${O}(n)$ models in two
  dimensions.
\newblock {\em Phys. Rev. Lett.}, 49:1062--1065, 1982.

\bibitem{Nie84}
B.~Nienhuis.
\newblock Coulomb gas description of {2D} critical behaviour.
\newblock {\em J. Statist. Phys.}, 34:731--761, 1984.

\bibitem{Orr47}
W.J.C. Orr.
\newblock Statistical treatment of polymer solutions at infinite dilution.
\newblock {\em Transactions of the Faraday Society}, 43:12--27, 1947.

\bibitem{St97}
J.~Michael Steele.
\newblock {\em Probability theory and combinatorial optimization}, volume~69 of
  {\em CBMS-NSF Regional Conference Series in Applied Mathematics}.
\newblock Society for Industrial and Applied Mathematics (SIAM), Philadelphia,
  PA, 1997.

\end{thebibliography}

\end{document}